\documentclass[11pt,a4paper]{article}
\usepackage[utf8]{inputenc}
\IfFileExists{lmodern.sty}{%
  \usepackage[T1]{fontenc}%
  \usepackage{lmodern}%
}{}
\IfFileExists{english.ldf}{\usepackage[english]{babel}}{}
\IfFileExists{microtype.sty}{\usepackage{microtype}}{}
\usepackage[margin=2.6cm]{geometry}
\usepackage{amsmath,amssymb,amsthm,mathtools}
\usepackage{enumitem}
\usepackage{booktabs,array}
\usepackage{graphicx}
\graphicspath{{figs/}{./}}
\newcommand{\figfile}[2][]{\includegraphics[#1]{#2}}
\usepackage{tikz}
\usetikzlibrary{decorations.pathreplacing}
\ifdefined\pdfsuppresswarningpagegroup\pdfsuppresswarningpagegroup=1\fi 
\usepackage{caption}
\usepackage{fancyhdr}
\usepackage[numbers,sort&compress]{natbib}
\usepackage[colorlinks=true,linkcolor=blue!45!black,citecolor=blue!45!black,urlcolor=blue!45!black]{hyperref}
\AtBeginDocument{\lefthyphenmin=2 \righthyphenmin=3}
\hypersetup{
  pdftitle={Persistence of a first-class particle in the two-species PushASEP},
  pdfauthor={Alikhan Aliaidar; Yuncherl Choi; Eunghyun Lee; Tamirlan Mukazhanov},
  pdfsubject={Capture, conditional joint positions and persistence for the two-species PushASEP},
  pdfkeywords={PushASEP, capture problem, persistence, first-passage, Bethe ansatz, determinantal formulas}
}

\numberwithin{equation}{section}
\newtheorem{theorem}{Theorem}[section]
\newtheorem{proposition}[theorem]{Proposition}
\newtheorem{lemma}[theorem]{Lemma}
\newtheorem{corollary}[theorem]{Corollary}
\theoremstyle{definition}
\newtheorem{definition}[theorem]{Definition}
\theoremstyle{remark}
\newtheorem{remark}[theorem]{Remark}

\newcommand{\Z}{\mathbb{Z}}
\newcommand{\C}{\mathbb{C}}
\newcommand{\PP}{\mathbb{P}}
\newcommand{\EE}{\mathbb{E}}
\newcommand{\cC}{\mathcal{C}}
\newcommand{\bA}{\mathbf{A}}
\newcommand{\bT}{\mathbf{T}}

\newcommand{\wa}{\mathsf{a}}
\newcommand{\wb}{\mathsf{b}}
\newcommand{\Inv}{\operatorname{Inv}}
\newcommand{\sgn}{\operatorname{sgn}}
\newcommand{\eps}{\varepsilon}
\newcommand{\FF}{\mathfrak{F}}
\newcommand{\GG}{\mathfrak{G}}
\newcommand{\setm}{\smallsetminus}
\newcommand{\diam}{\operatorname{diam}}

\title{\textbf{Persistence of a first-class particle\\
in the two-species PushASEP}}
\author{%
  Alikhan Aliaidar\textsuperscript{1}\quad
  Yuncherl Choi\textsuperscript{2}\quad
  Eunghyun Lee\textsuperscript{1,}\thanks{Corresponding author:
    \href{mailto:eunghyun.lee@nu.edu.kz}{\texttt{eunghyun.lee@nu.edu.kz}}}\quad
  Tamirlan Mukazhanov\textsuperscript{1}\\[4pt]
  {\small\textsuperscript{1}Nazarbayev University, Astana 010000, Kazakhstan}\\
  {\small\textsuperscript{2}Ingenium College, Kwangwoon University, Seoul 01891, Korea}%
}
\date{}

\begin{document}
\maketitle
\thispagestyle{fancy}

\begin{abstract}
\noindent A first-class particle starts to the right of $M=N-1$ second-class particles in the two-species PushASEP on $\Z$. Let $\tau$ be the first time this particle is no longer rightmost. For fixed $N$ and $p\in[0,1]$, we prove that $\PP_{Y(t)}(\tau>t)\sim\kappa_N\gamma_N(Y(t))t^{-M/2}$ whenever $y_N(t)-y_1(t)=o(\sqrt t)$. For fixed initial gaps the persistence exponent is $M/2$, and the mean capture time is finite exactly when $M\ge3$. The factor $\gamma_N$ is a Pitman--Stanley lattice-point count and a positive harmonic function for the killed process; $\kappa_N$ is a Gaussian integral independent of $p$. We give a joint local limit for the diffusively rescaled positions conditional on survival. When the second-class particles start at consecutive sites, an auxiliary non-collision process has the same survival probability and the same conditional distribution of the uppermost position at each time. For fixed initial sites the two conditional joint limits agree only for $N=2$. We also obtain crossing probabilities from $12\cdots2$ and $1\cdots12$ and tails of successive overtaking times.
\end{abstract}

{\small
\textbf{Keywords:} two-species PushASEP, persistence, capture time, Bethe ansatz, determinantal formulas, non-colliding walks.
\par}

\section{Introduction}

\subsection{The capture problem}\label{subsec:background}
How long does a diffusing prey stay ahead of $M$ predators that start on one side of it? For independent Brownian motions of equal diffusivity, including the prey, the survival exponent is $\beta_1=1/2$ for one predator and $\beta_2=3/4$ for two \cite{Bramson-Griffeath-1991,Krapivsky-Redner-1996,Redner-Krapivsky-1999}. Its growth for large $M$ is logarithmic: $\beta_M\sim\tfrac14\log M$ \cite[Theorem~3.3]{Li-Shao-2002}. For surveys of persistence see \cite{Bray-Majumdar-Schehr-2013,Aurzada-Simon-2015}.

Here the predators interact. A first-class particle starts to the right of $M$ second-class particles in the two-species PushASEP on $\Z$ \cite{Lee-2024}. We stop the process at the first change of this order. Until killing, its law is that of a free walk and an independent single-species PushASEP: a predator pushes a consecutive block on its right but is blocked by an occupied site on its left. In this representation killing occurs when the walk and the rightmost predator meet. In the two-species process it is an exchange, not occupation of the same site (Proposition~\ref{prop:killed-identification}).

For fixed initial gaps the survival exponent is $M/2$. Its origin is clearest when the predators start at consecutive sites. In that case the prey has the same probability of survival as a walk started above $M$ walks conditioned never to collide (Proposition~\ref{prop:packed-km}). For $N=M+1$ independent copies of the same nearest-neighbour walk started at distinct sites, non-collision contributes a factor of order $t^{-N(N-1)/4}$ \cite{Fisher,Karlin-McGregor-1959,Grabiner-1999,KOR}. Conditioning the lower $M$ walks contributes a Vandermonde weight of order $t^{M(M-1)/4}$ on the diffusive scale, leaving the exponent $M/2$. A monotone comparison gives the exponent for arbitrary fixed gaps, but not the prefactor.

For nonconsecutive initial sites, the auxiliary process generally has a different prefactor. The physical prefactor is nevertheless recovered by taking initial-coordinate differences of a Vandermonde-weighted auxiliary marginal (Remark~\ref{rem:initial-differences}). We identify it as a positive lattice-point count times a constant and obtain a joint local limit uniform over initial configurations of diameter $o(\sqrt t)$, in the sense of Theorem~\ref{thm:joint}. For $M\ge2$ the physical and auxiliary conditional joint limits differ even for fixed packed data.

The calculation uses the coordinate Bethe ansatz. The single-species TASEP kernel is a determinant \cite{Schutz-1997}; the ASEP kernel is a permutation sum of contour integrals \cite{Tracy-Widom-2008}. In multispecies systems the coefficients act on species words \cite{Tracy-Widom-2013,Lee-2020}. One-species counterparts include drop-push dynamics \cite{Schutz-Ramaswamy-Barma,Ali} and the PushASEP \cite{Borodin-Ferrari-2008}; for related determinantal diffusion models see \cite{Sasamoto-Wadati-1998}.

\subsection{The model}
Fix $N\ge2$. A state $(X,\pi)$ consists of ordered sites $x_1<\cdots<x_N$ in $\Z$ and a species word $\pi\in\{1,2\}^N$. An empty site has label $0$. Rank is counted from left to right.

\noindent\begin{minipage}{\linewidth}
\begin{definition}[Two-species PushASEP \cite{Lee-2024}]\label{def:model}
Species $2$ is first-class and species $1$ is second-class. Fix $p\in[0,1]$ and $q=1-p$. Each particle's attempts form an independent Poisson process of rate $1$. At an attempt by a particle of species $l$ at $x$, the update is as follows.
\begin{itemize}[leftmargin=2em]
\item[(R)] With probability $p$, the particle moves to the nearest site $z>x$ whose occupant has species $l'<l$. If $l'=1$ (so that $l=2$), the displaced species-$1$ particle moves to the nearest empty site to the right of $z$.
\item[(L)] With probability $q$, if the occupant of $x-1$ has species $l'<l$, the particle moves to $x-1$, and if $l'=1$ the species-$1$ particle that was at $x-1$ moves to $x$. If $l'\ge l$, nothing happens.
\end{itemize}
\end{definition}
\end{minipage}\par

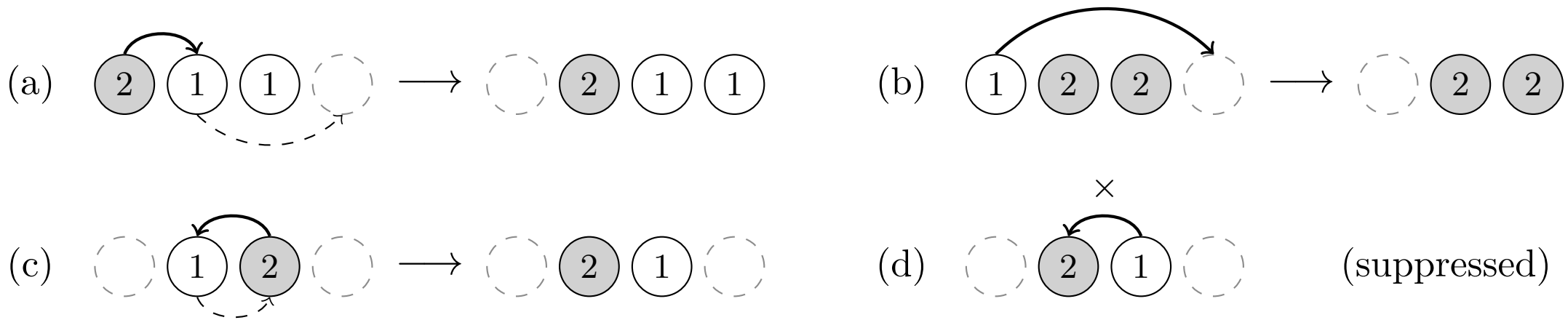
\begin{figure}[!htb]
\centering
\begin{tikzpicture}[x=6.6mm,y=7.2mm,
  site/.style={circle,draw,minimum size=5.4mm,inner sep=0pt,font=\footnotesize},
  two/.style={site,fill=black!18},
  one/.style={site,fill=white},
  emp/.style={site,dashed,draw=black!45}]
\node[font=\small] at (-1.3,0) {(a)};
\node[two] (a0) at (0,0) {2}; \node[one] (a1) at (1,0) {1}; \node[one] (a2) at (2,0) {1}; \node[emp] (a3) at (3,0) {};
\draw[->,thick] (a0.north) to[bend left=70] (a1.north);
\draw[->,dashed] (a1.south) to[bend right=45] (a3.south);
\node at (4.2,0) {$\longrightarrow$};
\node[emp] at (5.4,0) {}; \node[two] at (6.4,0) {2}; \node[one] at (7.4,0) {1}; \node[one] at (8.4,0) {1};
\node[font=\small] at (10.7,0) {(b)};
\node[one] (b0) at (12,0) {1}; \node[two] (b1) at (13,0) {2}; \node[two] (b2) at (14,0) {2}; \node[emp] (b3) at (15,0) {};
\draw[->,thick] (b0.north) to[bend left=45] (b3.north);
\node at (16.2,0) {$\longrightarrow$};
\node[emp] at (17.4,0) {}; \node[two] at (18.4,0) {2}; \node[two] at (19.4,0) {2}; \node[one] at (20.4,0) {1};
\node[font=\small] at (-1.3,-2.3) {(c)};
\node[emp] (c0) at (0,-2.3) {}; \node[one] (c1) at (1,-2.3) {1}; \node[two] (c2) at (2,-2.3) {2}; \node[emp] (c3) at (3,-2.3) {};
\draw[->,thick] (c2.north) to[bend right=70] (c1.north);
\draw[->,dashed] (c1.south) to[bend right=70] (c2.south);
\node at (4.2,-2.3) {$\longrightarrow$};
\node[emp] at (5.4,-2.3) {}; \node[two] at (6.4,-2.3) {2}; \node[one] at (7.4,-2.3) {1}; \node[emp] at (8.4,-2.3) {};
\node[font=\small] at (10.7,-2.3) {(d)};
\node[emp] (d0) at (12,-2.3) {}; \node[two] (d1) at (13,-2.3) {2}; \node[one] (d2) at (14,-2.3) {1}; \node[emp] (d3) at (15,-2.3) {};
\draw[->,thick] (d2.north) to[bend right=70] node[midway,above=0pt,font=\footnotesize] {$\times$} (d1.north);
\node[font=\small] at (18.2,-2.3) {(suppressed)};
\end{tikzpicture}
\caption{Right jumps: (a) a $2$ displaces a $1$, which moves to the first vacancy; (b) a $1$ jumps past two $2$'s to the first vacancy. Left jumps: (c) a $2$ exchanges with a $1$; (d) a $1$ is blocked.}
\label{fig:model}
\end{figure}

Right jumps may pass several particles; left jumps are nearest-neighbour (Figure~\ref{fig:model}). With a single species these rules give the PushASEP of \cite{Borodin-Ferrari-2008}, up to reflection. At $p=1$ only right jumps remain; for related priority rules on a ring see \cite{Bukh-Cox-2022,Ayyer-Martin-2023,AMW}.

\begin{remark}[Projections]\label{rem:proj}
Keeping only species $2$, or retaining all particles and erasing their labels, gives a single-species PushASEP under the rules of Definition~\ref{def:model} \cite{Lee-2024}. For a basic coupling (cf.\ \cite{Liggett-1999}), use the same independent right and left site clocks of rates $p$ and $q$ for both projections, ignoring rings at vacant sites. The first occupation set is then contained in the second; their difference consists of the sites carrying species $1$. Every change of word is a sequence of exchanges $12\to21$, never the reverse.
\end{remark}

We write $P_{(Y,\nu)}(X,\pi;t)$ for the probability that the system is in state $(X,\pi)$ at time $t$ given that it is in state $(Y,\nu)$ at time $0$. Throughout, $Y=(y_1,\dots,y_N)$ with $y_1<\dots<y_N$. We use $\PP_Y$ and $\EE_Y$ for probability and expectation from $Y$, with the initial species word, when present, specified by context.

\subsection{Main results}\label{subsec:results}
For $1\le m\le N$ let
\[
\wa_m:=\underbrace{2\cdots2}_{m-1}\,1\,\underbrace{2\cdots2}_{N-m},\qquad\qquad
\wb_m:=\underbrace{1\cdots1}_{m-1}\,2\,\underbrace{1\cdots1}_{N-m}
\]
denote the words with a single second-class particle or a single first-class particle, respectively, in position $m$. For the no-crossing and total-crossing marginals in Section~\ref{sec:marg}, the initial words are $\wa_1$ and $\wb_N$, respectively; the opposite orders $\wa_N$ and $\wb_1$ cannot change (Remark~\ref{rem:trivial}). Write $w_t$ for the law of a walk started at the origin, jumping right at rate $p$ and left at rate $q=1-p$.

\begin{figure}[!htb]
\centering
\begin{tikzpicture}[x=9mm,y=9mm,
  site/.style={circle,draw,minimum size=6.2mm,inner sep=0pt,font=\footnotesize},
  two/.style={site,fill=black!18}, one/.style={site,fill=white}, emp/.style={site,dashed,draw=black!40}]
\foreach \x in {-1,2,3,5,6,8,9} \node[emp] at (\x,0) {};
\node[one] (p1) at (0,0) {1}; \node[one] (p2) at (1,0) {1}; \node[one] (p3) at (4,0) {1}; \node[two] (w) at (7,0) {2};
\node[font=\footnotesize] at (0,-0.75) {$y_1$}; \node[font=\footnotesize] at (1,-0.75) {$y_2$};
\node[font=\footnotesize] at (4,-0.75) {$y_3$}; \node[font=\footnotesize] at (7,-0.75) {$y_4$};
\draw[->,thick] (p1.north) to[bend left=65] node[above,font=\footnotesize] {$p$} (p2.north);
\draw[->,dashed] (p2.north) to[bend left=65] (2,0.36);
\draw[->,thick] (p3.north) to[bend right=65] node[above,font=\footnotesize] {$q$} (3,0.36);
\draw[->,thick] (p3.north) to[bend left=65] node[above,font=\footnotesize] {$p$} (5,0.36);
\draw[->,thick] (w.north) to[bend right=65] node[above,font=\footnotesize] {$q$} (6,0.36);
\draw[->,thick] (w.north) to[bend left=65] node[above,font=\footnotesize] {$p$} (8,0.36);
\draw[decorate,decoration={brace,amplitude=4pt,mirror}] (-0.4,-1.15) -- (4.4,-1.15) node[midway,below=5pt,font=\footnotesize] {PushASEP with $N-1$ particles};
\node[font=\footnotesize] at (7,-1.55) {walk $W$};
\end{tikzpicture}
\caption{The killed representation for $N=4$. The grey walk starts at $y_N$; the white particles push to the right and block to the left. The first meeting with the walk is the killing time $\tau$.}
\label{fig:persistence}
\end{figure}
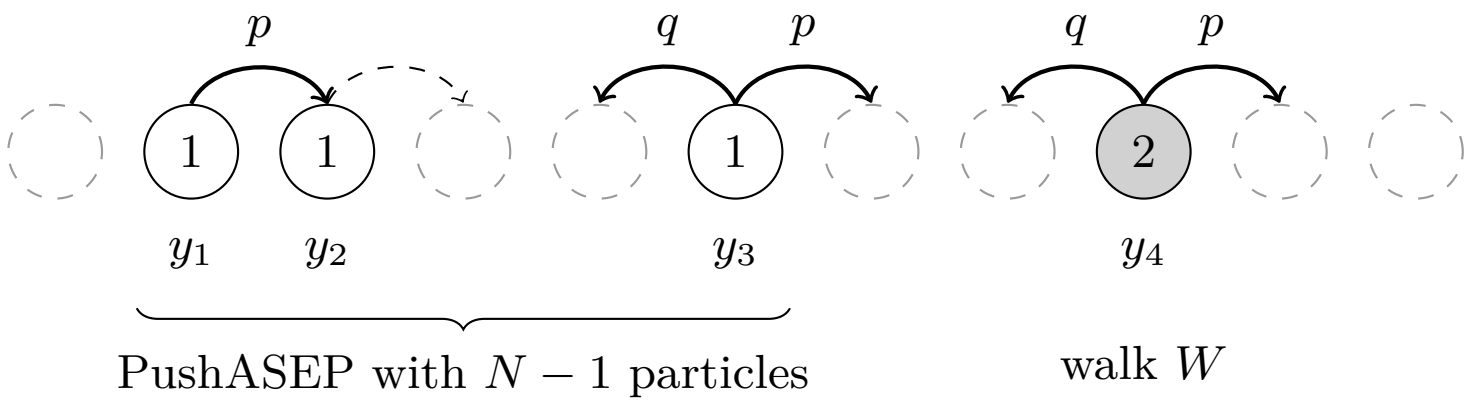

Let $W$ start at $y_N$ and let $\tau$ be its first meeting with the rightmost particle of an independent $(N-1)$-particle PushASEP started from $(y_1,\ldots,y_{N-1})$, where $y_{N-1}<y_N$. Proposition~\ref{prop:killed-identification} identifies this time in law with the first overtaking from $\wb_N$. For fixed $N$ and $p\in[0,1]$,
\begin{equation}\label{eq:intro-decay}
 \PP_Y(\tau>t)\sim\kappa_N\gamma_N(Y)t^{-(N-1)/2},\qquad y_N-y_1=o(\sqrt t).
\end{equation}
Equation~\eqref{eq:intro-decay} is Theorem~\ref{thm:decay}(c) for fixed $Y$ and Theorem~\ref{thm:joint} for moving initial data. Here $\gamma_N$ counts the integer chains in \eqref{eq:chains}, and $\kappa_N$ is the integral in \eqref{eq:kappa}. For example, at $Y=(0,2,5)$ the auxiliary gap factor is $15/2$, whereas $\gamma_3(Y)=9$ (Remark~\ref{rem:packed-general}). The same $\gamma_N$ gives the Doob transform for long survival (Corollaries~\ref{cor:harmonic} and~\ref{cor:doob}). For fixed initial sites the physical and auxiliary joint limits share their last marginal but differ for $N\ge3$ (Section~\ref{subsec:joint}). Corollary~\ref{cor:moments} gives the capture-time moments. Theorem~\ref{thm:overtaking} gives the individual overtaking-time tails also in the presence of several first-class particles; Corollary~\ref{cor:rank} gives the rank asymptotics for a single first-class particle.

Theorems~\ref{thm:redA} and~\ref{thm:redB} in Section~\ref{sec:marg} reduce the no-crossing and total-crossing endpoint sums from $\wa_1$ and $\wb_N$, respectively, to $k$-particle kernels.

\subsection{Relation to earlier work}\label{subsec:earlier}
The transition kernel is due to \cite{Lee-2024}. Its right and left updates impose the same collision boundary condition \cite[Remark~1]{Lee-2024}. The scattering operator is the specialization $(\widehat p,\widehat q)=(0,1)$ of \cite[equation~(1.5)]{Lee-Raimbekov-2022}, in that paper's rate convention. Reflection identifies $p=0$ with right-moving two-species TASEP; for $p>0$ only the dispersion changes, from $\xi-1$ to $p/\xi+(1-p)\xi-1$.
\begin{itemize}[leftmargin=2em,itemsep=1pt,topsep=2pt]
\item \emph{Theorem~\ref{thm:sN}.} At $p=0$, reflection gives the formula of \cite[Theorem~1.2]{Lee-2018} (Remark~\ref{rem:lee18}). This is also the case $k=1$ of \cite[Theorem~1.1]{Lee-2017}. The case $k=N-1$ of that theorem concerns the order $\wa_1$, but with the first-class particles at consecutive final sites, so it is a different quantity from $\FF_N$ in \eqref{eq:FN}.
\item \emph{Theorem~\ref{thm:matrix}.} For the corresponding ASEP words see \cite[Theorem~2.5 and Corollary~2.10]{Lee-Raimbekov-2022}. Remark~1.9 there omits the corresponding explicit lists of factors. Here each nonzero entry from $\wa_1$ or $\wb_N$ is a signed product given for every final rank.
\item \emph{Theorem~\ref{thm:nocross}.} At $p=0$, for any binary word, the reflected diagonal kernel agrees with \cite[Section~III, equation~(28)]{Chatterjee-Schutz-2010} and is the binary specialization of \cite[Theorem~3.1]{Lee-2020} (Remark~\ref{rem:nocross}). For the two-sided model such a determinant was anticipated in \cite[Section~5]{Lee-2024}.
\item \emph{Theorem~\ref{thm:totalcross}.} At $p=0$, the formulas agree with the cases of \cite[Proposition~2.1]{deGier-Mead-Wheeler-2023} in which one of the two types has a single particle (Remark~\ref{rem:dgmw}).
\end{itemize}
These are finite-time identities. Our asymptotics keep $N$ fixed, including at the TASEP endpoint. For packed second-class particles the equality of the surviving upper marginal in Proposition~\ref{prop:packed-km} also follows from interlacing constructions \cite{Borodin-Ferrari-2014,Assiotis-2018}, as explained in Section~\ref{subsec:packed}. The lattice-count polynomial is classical \cite{Stanley-Pitman-2002}. Here it is a positive harmonic function for the killed process and gives the normalization in the joint local limit for general initial gaps.

Section~\ref{sec:prelim} introduces the kernels, and Section~\ref{sec:decay} gives the persistence and overtaking results. Section~\ref{sec:trans} derives the species-word transition formulas; Section~\ref{sec:marg} treats the crossing marginals and their projections. The persistence proofs and open questions are in Sections~\ref{sec:persist-proofs} and~\ref{sec:open}. Appendices~\ref{app:matrix} and~\ref{app:lemmas} give the matrix-element and kernel proofs; Appendix~\ref{app:num} describes the numerical checks.

\section{Preliminaries}\label{sec:prelim}

\subsection{The Bethe-ansatz formula}
Let $e_\pi$ denote the basis vector indexed by a binary word $\pi$, and $[A]_{\pi,\nu}$ the coefficient of $e_\pi$ in $Ae_\nu$. Write $h(\xi)=1-\xi^{-1}$, $h_j=h(\xi_j)$, and, for $\alpha\ne\beta$, set
\begin{equation}\label{eq:factor}
 S_{\beta\alpha}=-h_\alpha/h_\beta,\qquad T_{\beta\alpha}=1-h_\alpha/h_\beta=1+S_{\beta\alpha}.
\end{equation}
For $1\le i<N$, let $s_i$ interchange slots $i,i+1$. The two-species scattering operator of \cite[Definition~6, equation~(57)]{Lee-2024} acts on these slots by
\begin{equation}\label{eq:Taction}
\bT_i(\beta,\alpha)\,e_\pi=\begin{cases}
-e_\pi+T_{\beta\alpha}\,e_{s_i\pi},&\text{if }\pi_i\pi_{i+1}=12,\\[2pt]
S_{\beta\alpha}\,e_\pi,&\text{otherwise.}
\end{cases}
\end{equation}
Identify $\sigma$ with the word $\sigma(1)\cdots\sigma(N)$: slot $i$ carries $\xi_{\sigma(i)}$ in \eqref{eq:Lee}. Starting from the identity, reach $\sigma$ by swaps in which slots $i,i+1$ contain $\alpha,\beta$ with $\alpha<\beta$ immediately before the swap. If the successive slots and pairs are $i_r,(\beta_r,\alpha_r)$, $1\le r\le\ell$, set
\begin{equation}\label{eq:Asigma}
\bA_\sigma:=\bT_{i_\ell}(\beta_\ell,\alpha_\ell)\cdots\bT_{i_2}(\beta_2,\alpha_2)\,\bT_{i_1}(\beta_1,\alpha_1).
\end{equation}
The consistency relations in \cite[Section~4]{Lee-2024} make this product independent of the reduced chain. We use $\bA_{\mathrm{id}}=I$.

\begin{theorem}[{\cite[Theorem~1, equation~(80)]{Lee-2024}}]\label{thm:Lee}
Let $0\le p\le1$, $q=1-p$, and let $\cC$ be a positively oriented circle $|\xi|=R>1$, $\oint_{\cC}=(2\pi i)^{-1}\int_{\cC}$ and $\eps(\xi)=p/\xi+q\xi-1$. For ordered $X,Y$, binary words $\pi,\nu$ of length $N$ and $t\ge0$,
\begin{equation}\label{eq:Lee}
P_{(Y,\nu)}(X,\pi;t)=\oint_{\cC}\!\cdots\oint_{\cC}\ \sum_{\sigma\in S_N}[\bA_\sigma]_{\pi,\nu}\ \prod_{i=1}^N\xi_{\sigma(i)}^{\,x_i-y_{\sigma(i)}-1}\ \prod_{i=1}^N e^{\eps(\xi_i)t}\ d\xi_1\cdots d\xi_N .
\end{equation}
\end{theorem}

In \cite[Definition~1 and Theorem~1]{Lee-2024}, $p$ and $q=1-p$ are jump probabilities, the initial state is arbitrary and the radii satisfy $1<R_1<\cdots<R_N$; the theorem imposes no separate strict inequality on $p$. It suffices here to invoke it for $0<p<1$. By \eqref{eq:factor}--\eqref{eq:Asigma}, every entry of $\bA_\sigma$ is a Laurent polynomial in $h_1,\ldots,h_N$. Each monomial in \eqref{eq:Lee} therefore gives a product of one-variable integrals of the form \eqref{eq:Fdef}, with no mixed poles. Each radius can be changed independently in $(1,\infty)$. On fixed contours the integrand is jointly continuous in $p$ and the contour variables. On the probability side, conditioning on the Poisson$(Nt)$ number of attempts gives a uniformly convergent mixture of polynomials in $p$: each conditional probability is bounded by one. Taking $p\downarrow0$ or $p\uparrow1$ proves the identity at both endpoints.

\subsection{The functions \texorpdfstring{$F_n$}{Fn}}
For integers $n,x$ define
\begin{equation}\label{eq:Fdef}
 F_n(x;t)=\oint_{\cC} h(\xi)^n \xi^{x-1} e^{\eps(\xi)t}\,d\xi .
\end{equation}
Below $W_t$ in the one-particle identities denotes a walk started at $0$. We suppress $t$ in $F_n$ when it is fixed.

\begin{lemma}\label{lem:F}
Let $n,x,a\in\Z$ and $t\ge0$.
\begin{enumerate}[label=\textup{(\alph*)},leftmargin=2.2em]
\item $F_0(x;t)=w_t(x):=\PP(W_t=x)=e^{-t}\sum_{s\ge\max(0,-x)}\dfrac{(pt)^{x+s}(qt)^{s}}{(x+s)!\,s!}$. At $p=1$, $w_t(x)=e^{-t}t^x/x!$ for $x\ge0$ and is zero otherwise; at $p=0$, $w_t(x)=e^{-t}t^{-x}/(-x)!$ for $x\le0$ and is zero otherwise, with the usual convention $0^0=1$.
\item $F_{n+1}(x;t)=F_n(x;t)-F_n(x-1;t)$.
\item $F_n(x;t)=\sum_{k\ge0}F_{n+1}(x-k;t)$, the series converging absolutely. In particular
\[
F_{-1}(x;t)=\PP(W_t\le x),\qquad F_{-2}(x;t)=\EE\big[(x+1-W_t)^+\big].
\]
\item $|F_n(x;t)|\le M_nR^{x}$ with $M_n=M_n(p,t,R):=\max_{\cC}|h^ne^{\eps t}|$; for $n\ge0$, $F_n(x;t)=\sum_{k=0}^n(-1)^k\binom nk w_t(x-k)$ and $\sum_{x\in\Z}|F_n(x;t)|\le2^n$.
\item $\displaystyle\sum_{x\le a}F_n(x;t)=F_{n-1}(a;t)$ for all $n$;\quad $\displaystyle\sum_{x>a}F_n(x;t)=-F_{n-1}(a;t)$ for $n\ge1$;\quad $\displaystyle\sum_{x>a}F_0(x;t)=1-F_{-1}(a;t)$.
\end{enumerate}
\end{lemma}

\begin{proof}
(a) On $\cC$ we have $e^{\eps(\xi)t}=e^{-t}\sum_{r,s\ge0}\frac{(pt)^r(qt)^s}{r!\,s!}\xi^{\,s-r}$, uniformly, and $\oint_\cC\xi^{x-1+s-r}d\xi=\delta_{r-s,x}$. The result is the law of the difference of two independent Poisson variables with parameters $pt$ and $qt$, which is the law of $W_t$.
(b) follows from $h(\xi)\xi^{x-1}=\xi^{x-1}-\xi^{x-2}$. The length of $\cC$ gives $|F_n(x;t)|\le M_nR^x$ directly from \eqref{eq:Fdef}.
(c) On $\cC$, $h(\xi)^{-1}=\sum_{k\ge0}\xi^{-k}$ uniformly. Termwise integration gives the identity, with absolute convergence since $\sum_{k\ge0}M_{n+1}R^{x-k}<\infty$. The two special cases follow from (a) and $\sum_{k\ge0}\mathbf 1_{\{W_t\le x-k\}}=(x+1-W_t)^+$.
(d) The remaining assertions follow from the binomial expansion in (b) and $\sum_xw_t(x)=1$.
(e) The first identity is (c). If $n\ge1$, then $F_{n-1}$ is absolutely summable by (d), hence $\sum_{x\in\Z}F_n(x)=\sum_{x\in\Z}[F_{n-1}(x)-F_{n-1}(x-1)]=0$, and the second identity follows from the first. The third one is (a) and (c).
\end{proof}

\subsection{Product form of the scattering factors}
For a word $\rho=\rho_1\cdots\rho_k$ of distinct letters in $\{1,\dots,N\}$ we set $\Inv(\rho):=\{(\rho_i,\rho_j): i<j,\ \rho_i>\rho_j\}$ and
\[
S_\rho:=\prod_{(\beta,\alpha)\in\Inv(\rho)}S_{\beta\alpha}.
\]
We denote by $\sgn(\rho)=(-1)^{|\Inv(\rho)|}$ the sign of $\rho$, and by $\operatorname{rk}_\rho(v)$ the rank of the letter $v$ among the letters of $\rho$ (the smallest letter has rank $1$).

\begin{lemma}\label{lem:Sprod}
For every word $\rho=\rho_1\cdots\rho_k$ of distinct letters,
\[
S_\rho=\sgn(\rho)\prod_{i=1}^k h_{\rho_i}^{\,i-\operatorname{rk}_\rho(\rho_i)}.
\]
In particular $S_\sigma=\sgn(\sigma)\prod_{i=1}^Nh_{\sigma(i)}^{\,i-\sigma(i)}$ for $\sigma\in S_N$.
\end{lemma}
\begin{proof}
By \eqref{eq:factor} each inversion contributes a sign, a factor $h_\alpha$ and a factor $h_\beta^{-1}$. Hence the exponent of $h_v$, $v=\rho_i$, is the number $L$ of letters larger than $v$ to the left of $v$ minus the number $s$ of letters smaller than $v$ to the right of $v$. If $s'$ is the number of letters smaller than $v$ to the left of $v$, then $L+s'=i-1$ and $s+s'=\operatorname{rk}_\rho(v)-1$, so $L-s=i-\operatorname{rk}_\rho(v)$.
\end{proof}

\begin{remark}[Single species]\label{rem:single}
If all particles have the same species, \eqref{eq:Taction} gives $[\bA_\sigma]_{\nu,\nu}=S_\sigma$, and Theorem~\ref{thm:Lee} together with Lemma~\ref{lem:Sprod} gives the determinant
\begin{equation}\label{eq:single}
P_Y(X;t)=\det\big[F_{i-j}(x_i-y_j;t)\big]_{i,j=1}^N
\end{equation}
of the PushASEP; cf.\ \cite{Borodin-Ferrari-2008}, and \cite[Section~2]{Ali} for $p=1$.
\end{remark}

\section{Persistence of the first-class particle: results}\label{sec:decay}
The proofs are in Section~\ref{sec:persist-proofs}. They use the diagonal transition kernel from Section~\ref{sec:trans} and, for successive overtakings, the projection in Section~\ref{sec:prob}.

\subsection{The problem}\label{subsec:problem}
The problem can be formulated without reference to species. Let $N\ge2$, let $W$ be a copy of the walk started at $y_N$, and let $x_1(t)<\dots<x_{N-1}(t)$ be an independent single-species PushASEP started from $y_1<\dots<y_{N-1}$, where $y_{N-1}<y_N$. Let
\[
\tau:=\inf\{t\ge0:\ x_{N-1}(t)=W(t)\}
\]
be the first meeting time (Figure~\ref{fig:persistence}). The independent system enters a cemetery state at this meeting; the two-species system is killed at its first overtaking, with its particles still at distinct sites.

\begin{proposition}[Identification of the killed process]\label{prop:killed-identification}
Start the two-species system from $(Y,\wb_N)$ and kill it at its first change of species order. Its surviving ordered positions have the same law as an independent single-species $(N-1)$-particle PushASEP and a walk $W$ started at $y_N$, killed at their first meeting. This holds for every $p\in[0,1]$. At an ordered state $Y$ the killing rate is
\[
k(Y)=\mathbf1_{\{y_N=y_{N-1}+1\}}\big(q+p\,\ell(Y)\big),
\]
where $\ell(Y)$ is the length of the consecutive block of second-class particles ending at $y_{N-1}$.
\end{proposition}

\begin{proof}
We compare the surviving transitions and the transitions to the cemetery state (the absorbing state which represents killing), clock by clock. The walk's right jump is always free in both descriptions. Its left jump is free unless $y_N-1=y_{N-1}$; in that case it meets the independent system, whereas the first-class particle exchanges order with its left neighbour in the two-species system. Both descriptions are killed at rate $q$.

A second-class particle's left jump is blocked only by another second-class particle, and therefore agrees with the single-species dynamics. For its right jump, let $b$ be the nearest vacant site when only the second-class particles are retained. If $b<y_N$, both systems have the same surviving update. If $b=y_N$, the independent system reaches the walk and is killed; in the two-species system the jumping particle passes the first-class particle and lands at $y_N+1$, changing the species order. The latter case occurs precisely for the $\ell(Y)$ clocks in the consecutive block immediately to the left of $y_N$, when that block is adjacent to the walk. Its total killing rate is $p\ell(Y)$. A right jump cannot pass the walk without meeting it in the independent description, since the walk's site is vacant for the single-species system.

Thus every surviving transition rate and the total killing rate coincide. The total attempt rate is $N$, so both killed chains are non-explosive and the common generator determines their law.
\end{proof}

Species exchanges cannot be reversed, so $\{\tau>t\}$ is the event that the species order at time $t$ is still $\wb_N$, and on this event the first-class particle is at $W(t)$. We set
\begin{align*}
s_N(x;t)&:=\PP\big(\tau>t,\ W(t)=x\big)=\sum_{x_1<\dots<x_{N-1}<x}P_{(Y,\wb_N)}\big((x_1,\dots,x_{N-1},x),\wb_N;t\big),\\
\mathcal S_N(t;Y)&:=\PP(\tau>t)=\sum_{x\in\Z}s_N(x;t).
\end{align*}
If all $N$ initial sites are consecutive, Proposition~\ref{prop:killed-identification} gives
$\mathcal S_N(t;Y)=1-[q+Mp]t+O(t^2)$ as $t\downarrow0$: the prey's left clock and each predator's right clock can kill on the first attempt. Thus for $M\ge2$ the finite-time survival depends on $p$, although the leading long-time term does not.

\subsection{Exact formula and asymptotics}\label{subsec:persist-results}
Write $e^{\wb}_j=j$ for $j<N$ and $e^{\wb}_N=N-1$.

\begin{theorem}\label{thm:sN}
For all $N\ge2$, $x\in\Z$ and $t\ge0$,
\begin{equation}\label{eq:sN}
s_N(x;t)=\det\big[F_{i-1-e^{\wb}_j}(x-y_j;t)\big]_{i,j=1}^N .
\end{equation}
\end{theorem}

\begin{remark}[The case $p=0$]\label{rem:lee18}
At $p=0$, Theorem~\ref{thm:sN} is, after the reflection $x\mapsto-x$, the formula of \cite[Theorem~1.2]{Lee-2018} for the right-moving TASEP started with a first-class particle to the left of $N-1$ second-class particles; it also agrees with the case $k=1$ of \cite[Theorem~1.1]{Lee-2017}. Write $\bar j=N+1-j$. In the reflected formula take the initial sites $-y_{\bar j}$, the observed site $-x$, and $\zeta_{\bar j}=\xi_j^{-1}$. The rational weight is the Vandermonde $\prod_{i<j}(\zeta_j-\zeta_i)$ times $\prod_j(1-\zeta_j)^{-(N+1-j-\delta_{j1})}$. The latter product becomes $\prod_jh_j^{-e^{\wb}_j}$. The Vandermonde becomes $\det[h_j^{\,i-1}]_{i,j=1}^N$: replacing $\zeta$ by $1-\zeta$ and reversing the columns each give $(-1)^{N(N-1)/2}$. Inverting each contour reverses its orientation and cancels the sign of its Jacobian, leaving $\xi_j^{\,x-y_j-1}e^{t(\xi_j-1)}d\xi_j$. Integrating the columns gives \eqref{eq:sN}.
\end{remark}

Let $\varphi(u)=(2\pi)^{-1/2}e^{-u^2/2}$ and define
\[
\Phi_0:=\varphi,\qquad \Phi_{-k}:=\varphi^{(k)}\ \ (k\ge1),\qquad \Phi_m(u):=\int_{-\infty}^u\Phi_{m-1}(r)\,dr\ \ (m\ge1),
\]
so that $\Phi_m'=\Phi_{m-1}$ for all $m\in\Z$ and $\Phi_1$ is the standard normal distribution function. Let
\begin{equation}\label{eq:Wkappa}
\omega_N(u):=\det\big[\Phi_{j-i}(u)\big]_{i,j=1}^N,\qquad \kappa_N:=\int_{\mathbb R}\omega_N(u)\,du ,
\end{equation}
and let $\gamma_N(Y)$ be the number of $(k_1,\dots,k_{N-1})\in\Z^{N-1}$ such that
\begin{equation}\label{eq:chains}
y_{N-j}\le k_j<k_{j-1}\quad(1\le j\le N-1),\qquad k_0:=y_N .
\end{equation}
Equivalently (put $m_j=y_N-k_j$), $\gamma_N(Y)$ is the number of integer sequences $0<m_1<\dots<m_{N-1}$ with $m_j\le y_N-y_{N-j}$ for all $j$. For example
\[
\gamma_2(Y)=y_2-y_1,\qquad \gamma_3(Y)=\binom{y_3-y_1}{2}-\binom{y_2-y_1}{2},
\]
and $\gamma_N(Y)=\binom{y_N-y_1}{N-1}$ when $y_j=y_1+j-1$ for $j<N$. The choice $k_j=y_{N-j}$ gives $\gamma_N(Y)>0$. Summing first over $k_1$ gives, with $\gamma_1=1$,
\begin{equation}\label{eq:gamma-recursion}
\gamma_N(y_1,\dots,y_N)=\sum_{k=y_{N-1}}^{y_N-1}
       \gamma_{N-1}(y_1,\dots,y_{N-2},k),
\end{equation}
so that $\gamma_N$ is a polynomial of degree $N-1$ in $Y$. Its combinatorial meaning is classical: set $M=N-1$, $d_j=y_{N-j+1}-y_{N-j}$, and let $\Pi_M(a)=\{z\in\mathbb R_{\ge0}^{M}:\sum_{i=1}^jz_i\le\sum_{i=1}^ja_i,\ 1\le j\le M\}$ be the Pitman--Stanley polytope \cite{Stanley-Pitman-2002}. With $m_0=0$, the bijection $z_j=m_j-m_{j-1}-1$ gives
\begin{equation*}
 \gamma_N(Y)=\#\big(\Pi_M(d-\mathbf1)\cap\Z^M\big).
\end{equation*}
The shift by $\mathbf1$ accounts for the strict inequalities in the integer chains.

All empty products are $1$. We write $\Delta_n(r):=\prod_{1\le k<l\le n}(r_l-r_k)$ for the Vandermonde determinant (so $\Delta_1=1$), and $v:=p-q$ for the drift of $W$.

\begin{theorem}\label{thm:decay}
Fix $N\ge2$, an ordered integer configuration $Y$, and $p\in[0,1]$, with $q=1-p$.
\begin{enumerate}[label=\textup{(\alph*)},leftmargin=2.2em]
\item For every $u\in\mathbb R$,
\begin{equation}\label{eq:Wformula}
\omega_N(u)=\frac{\varphi(u)}{\prod_{l=1}^{N-2}l!}\int_{r_1<\dots<r_{N-1}<u}\Delta_{N-1}(r)\,\Delta_N(r_1,\dots,r_{N-1},u)\prod_{k=1}^{N-1}\varphi(r_k)\,dr .
\end{equation}
In particular $\omega_N>0$, and
\begin{equation}\label{eq:kappa}
\kappa_N=\frac{1}{\prod_{l=1}^{N-2}l!}\int_{r_1<\dots<r_N}\Delta_{N-1}(r_1,\dots,r_{N-1})\,\Delta_N(r_1,\dots,r_N)\prod_{k=1}^{N}\varphi(r_k)\,dr\ \in(0,\infty).
\end{equation}
\item For every $A>0$,
\[
\lim_{t\to\infty}\ \sup_{x\in\Z:\ |x-vt|\le A\sqrt t}\Big|\,t^{N/2}s_N(x;t)-\gamma_N(Y)\,\omega_N\Big(\frac{x-vt}{\sqrt t}\Big)\Big|=0 .
\]
\item As $t\to\infty$,
\begin{equation}\label{eq:decay}
\mathcal S_N(t;Y)=\gamma_N(Y)\,\kappa_N\,t^{-(N-1)/2}\,\big(1+o(1)\big).
\end{equation}
\item For all $a<b$,
\[
\lim_{t\to\infty}\PP\Big(a<\frac{W(t)-vt}{\sqrt t}\le b\ \Big|\ \tau>t\Big)=\frac1{\kappa_N}\int_a^b\omega_N(u)\,du .
\]
\end{enumerate}
\end{theorem}

\begin{figure}[!htb]
\centering
\begin{minipage}[t]{.49\linewidth}\centering
\figfile[width=\linewidth]{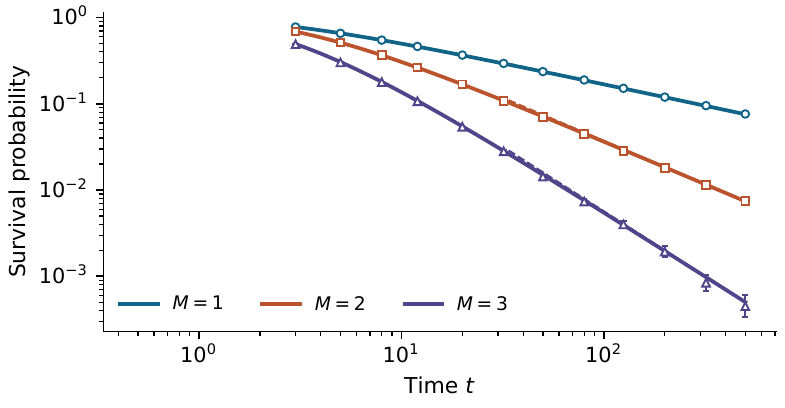}\\[-2pt]{\small(a) Survival probability}
\end{minipage}\hfill
\begin{minipage}[t]{.49\linewidth}\centering
\figfile[width=\linewidth]{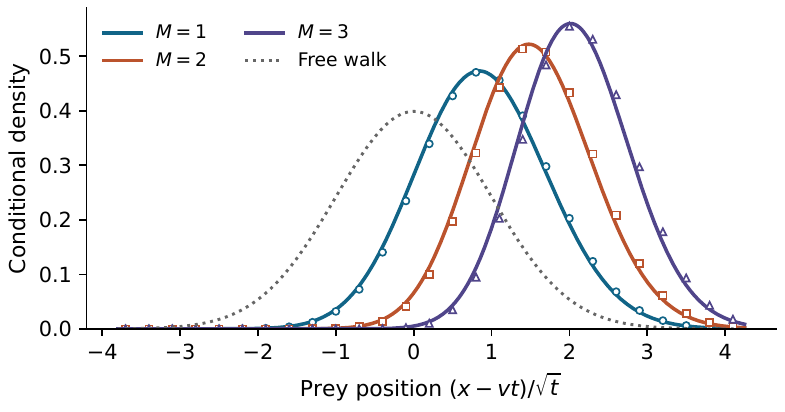}\\[-2pt]{\small(b) Conditional first-class density}
\end{minipage}
\caption{For $p=0.6$ and $Y=(0,3),(0,2,5),(0,2,3,6)$, respectively, $\gamma_N=3,9,16$.
(a) Solid: \eqref{eq:sN} summed over $x$; dashed: $\kappa_N\gamma_Nt^{-M/2}$, $M=N-1$. The leading term is independent of $p$; the finite-time survival need not be. Symbols: $10^5$ simulated paths per configuration, with pointwise 95\% Wilson intervals.
(b) Solid: $\omega_N/\kappa_N$; symbols: $\sqrt t\,s_N(x;t)/\mathcal S_N(t;Y)$ at $t=1600$; dotted: the standard normal density.}
\label{fig:survival}
\end{figure}

The free drivers have drift $v$ and variance rate $1$. After centring, $\omega_N/\kappa_N$ is independent of $p$. In \eqref{eq:Wformula} the integrated coordinates are auxiliary variables, not the positions of the second-class particles. Figure~\ref{fig:survival} displays the survival and the conditional marginal.

The lattice count is harmonic for the killed process.

\begin{corollary}\label{cor:harmonic}
The function $\gamma_N$ is harmonic for the process killed at time $\tau$: writing $Y(s):=(x_1(s),\dots,x_{N-1}(s),W(s))$, we have, for all $s\ge0$,
\[
\EE_Y\big[\gamma_N(Y(s));\ \tau>s\big]=\gamma_N(Y).
\]
\end{corollary}

Write $K_t(Y,X)=\PP_Y(Y(t)=X,\tau>t)$ for the killed kernel. For ordered $X,Y$, Proposition~\ref{prop:killed-identification} and Theorem~\ref{thm:nocross} give $K_t(Y,X)=\det[F_{e^{\wb}_i-e^{\wb}_j}(x_i-y_j;t)]$.

\begin{corollary}[Conditioning on long survival]\label{cor:doob}
The kernel
\begin{equation}\label{eq:doob-kernel}
K_t^\gamma(Y,X):=\frac{\gamma_N(X)}{\gamma_N(Y)}\,K_t(Y,X)
\end{equation}
is the transition kernel of a Markov process on the ordered configurations, and $\sum_{X}K^\gamma_t(Y,X)=1$. For each fixed $s\ge0$ and each bounded measurable functional $\Psi$ of the path up to time $s$,
\begin{equation}\label{eq:long-conditioning}
\lim_{\mathcal T\to\infty}\EE_Y[\Psi\mid\tau>\mathcal T]
=\frac{\EE_Y[\Psi\,\gamma_N(Y(s));\,\tau>s]}{\gamma_N(Y)},
\end{equation}
and the right-hand side is the expectation of $\Psi$ for the process with kernel $K^\gamma$ started from $Y$.
\end{corollary}

Here the path is observed up to a fixed time $s$, not at the conditioning horizon as in Theorem~\ref{thm:decay}(d).

\begin{remark}[The constants $\kappa_N$]\label{rem:kappa}
For independent standard normals $U,V$, \eqref{eq:kappa} gives $\kappa_2=\EE[(V-U)^+]=1/\sqrt\pi$. For $N=3$, set $a=r_2-r_1$, $b=r_3-r_2$. Integrating out the mean $(r_1+r_2+r_3)/3$ in \eqref{eq:kappa} and substituting $b=as$ gives
\[
 \kappa_3=\frac{1}{2\pi\sqrt3}\int_{a,b>0}a^2b(a+b)e^{-(a^2+ab+b^2)/3}\,da\,db
 =\frac{27}{2\pi\sqrt3}\int_0^\infty\frac{s(1+s)}{(1+s+s^2)^3}\,ds
 =\frac{3\sqrt3}{4\pi}.
\]
The last integrand is $-\tfrac16$ times the derivative of $(2s+1)/(1+s+s^2)^2$. Numerical integration gives $\kappa_4\approx0.35437$ and $\kappa_5\approx0.34010$.
\end{remark}

\subsection{Packed initial data and non-colliding walks}\label{subsec:packed}
Call the initial data \emph{packed} when the $M=N-1$ second-class particles start at consecutive sites. In the homogeneous interlacing array of \cite{Borodin-Ferrari-2014,Assiotis-2018}, the right edge pushes to the right and blocks to the left. With positive constant right and left rates and Gibbs initial data, the longest row evolves as the corresponding walks conditioned never to collide. If this row is packed, interlacing ($a_i\le b_i<a_{i+1}$ for successive rows $a,b$, with $a$ the longer row) fixes every shorter row. The edge and the longest row then share their rightmost trajectory; adding the same independent walk above them gives the equality of surviving upper marginals below. See also \cite[Remark~2.23]{Borodin-Ferrari-2014} for two-sided dynamics and \cite{OConnell-Yor-2002,KOR} for representations of non-colliding walks. Section~\ref{subsec:packed-proof} gives a direct determinant proof, including $p=0,1$.

\begin{proposition}[Reduction to non-colliding walks for packed initial data]\label{prop:packed-km}
Let $M=N-1$, $Y_-=(y_1,\dots,y_M)$, and suppose that $y_j=y_1+j-1$ for $1\le j\le M$, with $y_N>y_1+M-1$. Then, for every $x\in\Z$ and $t\ge0$,
\begin{equation}\label{eq:packed-km}
s_N(x;t)=\sum_{z_1<\cdots<z_M<x}
 \frac{\Delta_M(z_1,\dots,z_M)}{\Delta_M(Y_-)}
 \det\big[w_t(z_i-y_j)\big]_{i,j=1}^N,
 \qquad z_N=x.
\end{equation}
This is the joint mass of survival and the endpoint $x$ of an independent upper walk above $M$ non-colliding walks (the lower walks have the Doob transform by $\Delta_M$). The equality concerns the upper marginal, not the law of the lower positions. Moreover,
\begin{equation}\label{eq:packed-gamma}
\gamma_N(Y)=\binom{y_N-y_1}{M}
 =\frac{\Delta_N(Y)}{M!\,\Delta_M(Y_-)}.
\end{equation}
\end{proposition}

\begin{remark}[Packed and general initial data]\label{rem:packed-general}
For packed data, \eqref{eq:packed-km} recovers the marginal and survival asymptotics of Theorem~\ref{thm:decay} from the non-collision formula \cite{Karlin-McGregor-1959,KOR,Grabiner-1999}. For arbitrary fixed $Y$, the auxiliary walks on the right of \eqref{eq:packed-km} have survival prefactor $\kappa_N\Delta_N(Y)/(M!\,\Delta_M(Y_-))$ (Section~\ref{subsec:packed-proof}). For $N=3$, let $g=y_2-y_1$ and $d=y_3-y_2$. By \eqref{eq:chains},
\[
 \gamma_3(Y)-\frac{\Delta_3(Y)}{2\Delta_2(Y_-)}
 =\binom{g+d}{2}-\binom{g}{2}-\frac{d(g+d)}2
 =\frac{d(g-1)}2 .
\]
Thus the two gap factors agree exactly when the lower pair is packed. Taking finite differences across starting configurations gives the physical gap factor instead (Remark~\ref{rem:initial-differences}).

For $N=2$ both killed kernels describe two independent walks stopped at their first meeting, since $\Delta_1=1$. For $M\ge2$ they can differ even from packed data: take $y_N>y_M+1$. From packed $Y_-$, moving $y_1$ alone one step left has rate $q$ in the PushASEP and $Mq$ in the auxiliary process; moving $y_M$ alone one step right has rates $p$ and $Mp$. In either move the lower Vandermonde ratio is $\prod_{j=1}^{M-1}(j+1)/j=M$. At least one rate differs for every $p\in[0,1]$.
\end{remark}

A comparison with two packed configurations gives the survival exponent for arbitrary fixed gaps.
For $1\le i\le M$ set $l_i=y_1+i-1$ and $r_i=y_M-M+i$. Run the three \emph{unkilled} predator systems from $l\le Y_-\le r$ with common rank clocks and directions, and use the same independent prey. On a right attempt at rank $j$, increment $x_j$ and then set $x_i'=\max\{x_i,x_{i-1}'+1\}$ for $i>j$. A left attempt replaces $x_j$ by $\max\{x_j-1,x_{j-1}+1\}$ (by $x_1-1$ if $j=1$). All other coordinates are unchanged. Since increment, decrement and maximum are nondecreasing in their arguments, induction over ranks and then over clock rings preserves the three-way order. Each front moves by at most one per ring, and its clock rings and the prey's jumps are almost surely disjoint. Each gap $W-x_M$ therefore hits zero before becoming negative. The rightmost front meets the prey no later than the middle one, and the middle one no later than the leftmost one. Thus
\[
 \mathcal S_N(t;(r,y_N))\le \mathcal S_N(t;Y)\le\mathcal S_N(t;(l,y_N)).
\]
For fixed $Y$ the two endpoints are packed, so \eqref{eq:packed-km} and non-collision asymptotics already give $\mathcal S_N(t;Y)=\Theta(t^{-M/2})$. The comparison does not identify the amplitude between the two bounds or the joint conditional density.

\subsection{Joint positions at the survival horizon}\label{subsec:joint}
Formula \eqref{eq:Wformula} describes only the last coordinate. Keeping all rows of $K_t$ gives the following joint limit.

\begin{theorem}[Joint local limit for the surviving positions]\label{thm:joint}
Fix $N\ge2$ and $p\in[0,1]$, and write $M=N-1$, $e=e^{\wb}=(1,\ldots,M,M)$. Allow the deterministic initial configuration $Y=Y^{[t]}$ to vary with the observation horizon. For an ordered endpoint $X$ put $u_i=(x_i-y_1-vt)/\sqrt t$ and define on the closed chamber $u_1\le\cdots\le u_N$
\begin{equation}\label{eq:joint-profile}
 Q_N(u)=\det\big[\Phi_{j-e_i}(u_i)\big]_{\substack{1\le i\le N\\0\le j\le M}},
\end{equation}
where columns are ordered by increasing $j$. If the diameter $\diam(Y):=y_N-y_1$ of the initial configuration is $o(\sqrt t)$, then
\begin{align}
 \frac{t^{(2N-1)/2}}{\gamma_N(Y)}K_t(Y,X)&=Q_N(u)+o(1),\label{eq:joint-llt}\\
 \frac{t^{M/2}}{\gamma_N(Y)}\mathcal S_N(t;Y)&=\kappa_N+o(1).\label{eq:joint-survival}
\end{align}
For every $\eta(t)\downarrow0$ these estimates are uniform over $Y$ with $\diam(Y)\le\eta(t)\sqrt t$, and the first is uniform over endpoints with $u$ in any fixed compact set. In particular, conditional on $\tau>t$, the rescaled ordered positions have local masses
\[
 t^{N/2}\frac{K_t(Y,X)}{\mathcal S_N(t;Y)}=\frac{Q_N(u)}{\kappa_N}+o(1)
\]
on compact sets and converge weakly to the law with this density on the ordered chamber. The density of its last coordinate is $\omega_N/\kappa_N$.
\end{theorem}

For fixed initial sites, Section~\ref{subsec:packed-proof} identifies the auxiliary conditional joint limit: its density is proportional on the ordered chamber to $P(u)\prod_i\varphi(u_i)$, where $u_-=(u_1,\ldots,u_M)$ and $P(u)=\Delta_M(u_-)\Delta_N(u)$. Its last marginal is $\omega_N/\kappa_N$ by \eqref{eq:Wformula}--\eqref{eq:kappa}. For $N=2$ the joint limits coincide, since $Q_2(u)=(u_2-u_1)\varphi(u_1)\varphi(u_2)$. The two profiles are related algebraically. On $\mathbb R^N$ let $\mathcal I_i$ integrate the $i$th argument from $-\infty$ to $u_i$, holding the other arguments fixed. Then
\[
 Q_N(u)=\left(\prod_{i=1}^{N}\mathcal I_i^{\,M-e_i}\right)
        \left[\Delta_N(u)\prod_{i=1}^N\varphi(u_i)\right].
\]
Indeed, the bracket is $\det[\Phi_{j-M}(u_i)]_{1\le i\le N,\,0\le j\le M}$: the Gaussian derivative signs cancel the sign of reversing the columns. Applying the indicated primitives row by row gives \eqref{eq:joint-profile}. The input here is the alternating function on $\mathbb R^N$, not the auxiliary density restricted to the ordered chamber.

For $N\ge3$ the joint limits differ. Differentiating the Gaussian identity $\varphi'+u\varphi=0$, or integrating from $-\infty$ with vanishing boundary terms, gives $\Phi_m''+u\Phi_m'+(1-m)\Phi_m=0$ for every integer $m$. Put $\partial_i=\partial/\partial u_i$ and $\mathcal D=\sum_i(\partial_i^2+u_i\partial_i)$. In each determinant product in \eqref{eq:joint-profile}, the sum of the coefficients $1-m$ is $\sum_i(1+e_i)-\sum_{j=0}^M j=2N-1$. Thus $(\mathcal D+2N-1)Q_N=0$. On the auxiliary profile the same operator gives
\[
 (\mathcal D+2N-1)\Big(P\prod_i\varphi(u_i)\Big)
 =\Big(\sum_i\partial_i^2P+(M-M^2)P\Big)\prod_i\varphi(u_i),
\]
since $P$ is homogeneous of degree $M^2$. For $M\ge2$ the highest-degree term is nonzero, so the two profiles cannot be proportional on the open chamber.

For $N=3$ the difference is visible on the lower-particle contact face $u_1=u_2<u_3$ (the reflecting face):
\[
 Q_3(0,0,1)=\varphi(0)\varphi(1)\big(\varphi(0)+\tfrac12\big)>0 .
\]
The auxiliary density vanishes on that face. By continuity the densities differ also at interior points near $(0,0,1)$; neither law puts mass on the face.

\begin{figure}[!htb]
\centering
\figfile[width=0.69\linewidth]{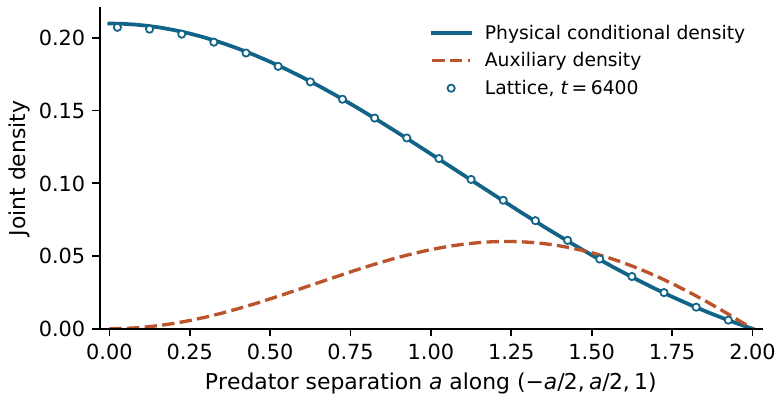}
\caption{A slice through the three-particle conditional limit densities at $u=(-a/2,a/2,1)$, $0<a<2$. Solid: $Q_3/\kappa_3$ for the PushASEP. Dashed: the auxiliary non-collision density, the normalized integrand of \eqref{eq:kappa} (Section~\ref{subsec:packed-proof}). Open circles: $t^{3/2}K_t(Y,X)/\mathcal S_3(t;Y)$ for $p=1/2$, $Y=(0,2,5)$ and $t=6400$. The last marginals agree, but the joint densities differ near $a=0$.}
\label{fig:joint}
\end{figure}

Figure~\ref{fig:joint} compares the two densities near contact. For $N=2$ and $y_2-y_1=\lfloor t^{1/4}\rfloor$, \eqref{eq:joint-survival} gives $\mathcal S_2(t;Y)\sim t^{-1/4}/\sqrt\pi$: the exponent $1/2$ for a fixed gap does not describe a growing gap. At the diffusive boundary $(y_2-y_1)/\sqrt t\to a>0$, Corollary~\ref{cor:cross}(a) and the central limit theorem for the difference walk (variance $2t$) give $\mathcal S_2(t;Y)\to\operatorname{erf}(a/2)$, not the linear approximation $a/\sqrt\pi$.

\subsection{Comparison with independent predators}\label{subsec:capture}
For an independent Brownian prey and $M$ independent Brownian predators initially on one side of it, all with the same diffusivity, the survival exponent $\beta_M$ is determined by a spherical eigenvalue problem \cite{Bramson-Griffeath-1991,Krapivsky-Redner-1996,Redner-Krapivsky-1999}. Table~\ref{tab:capture} compares this exponent with the fixed-gap lattice exponent $M/2$.

\begin{table}[!htb]
\centering
\begin{tabular}{lccccc}
\toprule
number of predators $M$ & $1$ & $2$ & $3$ & $4$ & $10$\\
\midrule
independent predators, $\beta_M$ & $1/2$ & $3/4$ & $0.91288$ & $1.03$ & $1.4$\\
PushASEP predators, $M/2$ & $1/2$ & $1$ & $3/2$ & $2$ & $5$\\
\bottomrule
\end{tabular}
\caption{Survival exponents for fixed initial gaps. For $M=3$ independent predators started together, the value is rounded from \cite[equation~(27)]{Biroli-2026}. The estimates for $M=4,10$ are due to \cite{Bramson-Griffeath-1991} and reported in \cite{Krapivsky-Redner-1996}. The PushASEP values are exact.}
\label{tab:capture}
\end{table}

\begin{corollary}[Moments of the capture time]\label{cor:moments}
Let $s>0$. For every initial configuration $Y$, $\EE_Y[\tau^s]<\infty$ if and only if $s<(N-1)/2$. In particular, the mean capture time $\EE_Y[\tau]$ is finite if and only if $N\ge4$, that is, if and only if the walk is pursued by at least three particles.
\end{corollary}

Indeed, $\EE_Y[\tau^s]=s\int_0^\infty t^{s-1}\PP_Y(\tau>t)\,dt$ and \eqref{eq:decay} has a strictly positive prefactor. For independent Brownian predators the mean is infinite for $M\le3$ \cite{Bramson-Griffeath-1991}; it is finite for $M\ge5$ by \cite{Li-Shao-2001} and for $M=4$ by \cite{Ratzkin-Treibergs-2009}. For distinct fixed starting positions on one side, using the same Brownian increments bounds the capture time between those for predators all started at the nearest or the farthest initial position. Brownian scaling leaves the moment threshold unchanged.

\subsection{Successive overtakings and several first-class particles}\label{subsec:overtaking}
Consider the two-species PushASEP started from an arbitrary state with at least one first-class particle. Let $Z(t)$ be the position of the leftmost first-class particle at time $t$, and let $R(t)$ be the number of second-class particles to the left of $Z(t)$. Equivalently, $R(t)$ is the length of the leading block of $1$'s in the species word at time $t$. Each elementary exchange $12\to21$ can only shorten this block or leave it unchanged (Remark~\ref{rem:proj}), regardless of which first-class particle occupies $Z(t)$. For $1\le m\le R(0)$ the stopping time
\[
\tau_m:=\inf\{t\ge0:\ R(t)<m\}
\]
is the first time the number of second-class particles strictly to the left of $Z(t)$ falls below $m$. Then $\tau_{R(0)}\le\dots\le\tau_2\le\tau_1$ and $\{\tau_m>t\}=\{R(t)\ge m\}$. For the initial state $(Y,\wb_N)$ we have $R(0)=N-1$ and $\tau_{N-1}=\tau$.

\begin{theorem}[Successive overtakings]\label{thm:overtaking}
Let $p\in[0,1]$ and $m\ge1$. Consider the two-species PushASEP started from a state with $R(0)\ge m$. Let $y_1'<\dots<y_m'$ be the initial positions of the $m$ leftmost second-class particles, let $z_0>y_m'$ be the initial position of the leftmost first-class particle, and put $Y^{(m)}:=(y_1',\dots,y_m',z_0)$.
\begin{enumerate}[label=\textup{(\alph*)},leftmargin=2em]
\item Send the vector of the $m$ leftmost second-class positions and the leftmost first-class position to a cemetery state at $\tau_m$. This killed process has the same law as an independent $m$-particle PushASEP from $(y_1',\ldots,y_m')$ and walk from $z_0$, sent to the same state at their first meeting.
\item For all $x\in\Z$ and $t\ge0$,
\[
\PP\big(R(t)\ge m,\ Z(t)=x\big)=\det\big[F_{i-1-e_j}\big(x-Y^{(m)}_j;t\big)\big]_{i,j=1}^{m+1},\qquad e=(1,2,\dots,m,m),
\]
which is the right-hand side of \eqref{eq:sN} with $N$ replaced by $m+1$ and $Y$ by $Y^{(m)}$. In particular $\PP(\tau_m>t)=\mathcal S_{m+1}\big(t;Y^{(m)}\big)$.
\item As $t\to\infty$,
\[
\PP(\tau_m>t)=\PP\big(R(t)\ge m\big)=\gamma_{m+1}\big(Y^{(m)}\big)\,\kappa_{m+1}\,t^{-m/2}\,\big(1+o(1)\big),
\]
and, given $\{R(t)\ge m\}$, the rescaled position $(Z(t)-(p-q)t)/\sqrt t$ converges in law to the distribution with density $\omega_{m+1}/\kappa_{m+1}$.
\end{enumerate}
\end{theorem}

Only the $m$ leftmost second-class positions and the leftmost first-class position enter this law. The identification ends at killing: from $((0,1),12)$ a right jump of the left particle gives the physical state $((1,2),21)$, whereas in the independent description $(x_1,W)=(1,1)$ at contact. The theorem specifies the separate laws of the $\tau_m$, not their dependence.
For a single first-class particle, differencing the nested tail probabilities gives the rank probabilities.

\begin{corollary}\label{cor:rank}
Start the system from $(Y,\wb_N)$. Then, as $t\to\infty$:
\begin{enumerate}[label=\textup{(\alph*)},leftmargin=2em]
\item for $1\le m\le N-1$, the species order at time $t$ is $\wb_{m+1}$ with probability
\[
\PP\big(R(t)=m\big)=\gamma_{m+1}(y_1,\dots,y_m,y_N)\,\kappa_{m+1}\,t^{-m/2}\,\big(1+o(1)\big);
\]
\item the number of second-class particles strictly to the left of the first-class particle has the mean
\[
\EE\big[R(t)\big]=\frac{y_N-y_1}{\sqrt{\pi t}}\,\big(1+o(1)\big);
\]
\item for $1\le m\le N-1$ and $s>0$, $\EE[\tau_m^s]<\infty$ if and only if $s<m/2$; in particular $\EE[\tau_m]<\infty$ if and only if $m\ge3$.
\end{enumerate}
\end{corollary}

For $m<N-1$, subtract the smaller-order probability $\PP(R(t)\ge m+1)$ from $\PP(R(t)\ge m)$; for $m=N-1$ the subtracted term is zero. Summing the tails gives (b), since $\gamma_2(y_1,y_N)\kappa_2=(y_N-y_1)/\sqrt\pi$. Part~(c) follows by tail integration.

\section{Method: matrix elements and transition probabilities}\label{sec:trans}
We evaluate the matrix entries for the initial orders $\wa_1$ and $\wb_N$, the diagonal kernel for any binary word, and the two total-crossing kernels.

\subsection{Matrix elements}\label{sec:matrix}
For $\sigma\in S_N$ and a letter $v$ we denote by $\sigma\setm v$ the word of length $N-1$ obtained from $\sigma(1)\cdots\sigma(N)$ by deleting the letter $v$.

\begin{theorem}\label{thm:matrix}
Let $\sigma\in S_N$ and $1\le m\le N$.
\begin{enumerate}[label=\textup{(\alph*)},leftmargin=2.2em]
\item Let $r:=\sigma^{-1}(1)$ be the position of the letter $1$ in $\sigma$. Then
\begin{equation}\label{eq:matA}
[\bA_\sigma]_{\wa_m,\wa_1}=\begin{cases}
\displaystyle S_{\sigma\setm1}\prod_{i=1}^{r-1}T_{\sigma(i)1},& m=r,\\[12pt]
\displaystyle -\,S_{\sigma\setm1}\prod_{i=1}^{m-1}T_{\sigma(i)1}\prod_{i=m+1}^{r-1}S_{\sigma(i)1},& m<r,\\[12pt]
0,& m>r.
\end{cases}
\end{equation}
\item Let $r:=\sigma^{-1}(N)$ be the position of the letter $N$ in $\sigma$. Then
\begin{equation}\label{eq:matB}
[\bA_\sigma]_{\wb_m,\wb_N}=\begin{cases}
\displaystyle S_{\sigma\setm N}\prod_{i=r+1}^{N}T_{N\sigma(i)},& m=r,\\[12pt]
\displaystyle -\,S_{\sigma\setm N}\prod_{i=m+1}^{N}T_{N\sigma(i)}\prod_{i=r+1}^{m-1}S_{N\sigma(i)},& m>r,\\[12pt]
0,& m<r.
\end{cases}
\end{equation}
\end{enumerate}
\end{theorem}

The proof consists in choosing a convenient chain of elementary steps for $\sigma$; it is given in Appendix~\ref{app:matrix}.

\begin{remark}[Special cases]\label{rem:special}
(i) For $m=1$, \eqref{eq:matA} reads: $[\bA_\sigma]_{\wa_1,\wa_1}=S_\sigma$ if $\sigma(1)=1$, and $[\bA_\sigma]_{\wa_1,\wa_1}=-\prod S_{\beta\alpha}$, the product being over $\Inv(\sigma)\setm\{(\sigma(1),1)\}$, if $\sigma(1)\neq1$. Indeed $\Inv(\sigma)$ is the disjoint union of $\Inv(\sigma\setm1)$ and $\{(\sigma(i),1): i<r\}$.

(ii) For $m=1$, \eqref{eq:matB} reads: $[\bA_\sigma]_{\wb_1,\wb_N}=\prod_{k=1}^{N-1}T_{Nk}\cdot S_{\sigma\setm N}$ if $\sigma(1)=N$, and $[\bA_\sigma]_{\wb_1,\wb_N}=0$ otherwise. Thus only $(N-1)!$ of the $N!$ terms of \eqref{eq:Lee} contribute to the total crossing probability $\wb_N\to\wb_1$, and similarly for $\wa_1\to\wa_N$.
\end{remark}

\begin{remark}[The orders $2\cdots21$ and $21\cdots1$]\label{rem:trivial}
If $\nu=\wa_N$ or $\nu=\wb_1$, no word reached from $e_\nu$ ever has the letters $12$ in adjacent slots, so that by \eqref{eq:Taction} $[\bA_\sigma]_{\pi,\nu}=\delta_{\pi,\nu}S_\sigma$ and $P_{(Y,\nu)}(X,\nu;t)$ is the single-species determinant \eqref{eq:single}, in agreement with Remark~\ref{rem:proj}.
\end{remark}

\subsection{No crossing}
\begin{theorem}\label{thm:nocross}
For any binary word $\nu$, define the integer exponent vector $e=(e_1,\ldots,e_N)$ by $e_1=0$ and $e_{i+1}-e_i=\mathbf1_{\{(\nu_i,\nu_{i+1})\ne(1,2)\}}$. For $p\in[0,1]$ and ordered $X,Y$ the joint mass of the endpoint and no exchange is
\begin{equation}\label{eq:anyword}
 P_{(Y,\nu)}(X,\nu;t)=\det\big[F_{e_i-e_j}(x_i-y_j;t)\big]_{i,j=1}^N .
\end{equation}
Only differences of $e$ occur. For the two extremal orders we may take $e^{\wa}_i=i+\delta_{i1}$ and $e^{\wb}_i=i-\delta_{iN}$, respectively, obtaining
\begin{align}
P_{(Y,\wa_1)}(X,\wa_1;t)&=\det\big[F_{e^{\wa}_i-e^{\wa}_j}(x_i-y_j;t)\big]_{i,j=1}^N,\label{eq:detA}\\
P_{(Y,\wb_N)}(X,\wb_N;t)&=\det\big[F_{e^{\wb}_i-e^{\wb}_j}(x_i-y_j;t)\big]_{i,j=1}^N.\label{eq:detB}
\end{align}
\end{theorem}

\begin{proof}
Every off-diagonal update $12\to21$ increases the number of inversions of the species word; no path using such an update can return to $\nu$. Along any elementary chain for $\sigma$, the diagonal coefficient is therefore the product of the diagonal factors for the unchanged word. Put $\delta_i=e_{i+1}-e_i$. If spectral letters $\alpha<\beta$ are exchanged at slots $i,i+1$, equations~\eqref{eq:factor} and~\eqref{eq:Taction} give the factor $-(h_\alpha/h_\beta)^{\delta_i}$. The ratio of
\[
 c_\nu(w)=\sgn(w)\prod_{i=1}^N h_{w(i)}^{\,e_i-e_{w(i)}}
\]
after and before this step is exactly this factor. Since $c_\nu(\mathrm{id})=1$, the diagonal coefficient is $c_\nu(\sigma)$. Substitution in \eqref{eq:Lee} gives a signed permutation sum of products of one-variable integrals $F_{e_i-e_{\sigma(i)}}(x_i-y_{\sigma(i)};t)$, the determinant \eqref{eq:anyword}. For a fixed species word, ending in that word is equivalent to no exchange. To obtain the conditional endpoint law one divides by the sum of \eqref{eq:anyword} over ordered endpoints.
\end{proof}

\begin{remark}\label{rem:nocross}
For $N=2$ both \eqref{eq:detA} and \eqref{eq:detB} reduce to the Karlin--McGregor determinant $w_t(x_1-y_1)w_t(x_2-y_2)-w_t(x_1-y_2)w_t(x_2-y_1)$ for two independent random walks \cite{Karlin-McGregor-1959}; see Section~\ref{sec:prob} for the explanation. Determinants of the form $\det[F_{e_i-e_j}(x_i-y_j)]$ thus interpolate between the Karlin--McGregor formula ($e$ constant) and the Sch\"utz-type formula \eqref{eq:single} ($e_i=i$); compare the kernels of \cite{Dieker-Warren-2008}. At $p=0$ the comparison applies to every binary word in \eqref{eq:anyword}. Put $\bar i=N+1-i$ and reflect $(x_i,y_i,\nu_i)$ to $(-x_i,-y_i,\nu_i)$ in slot $\bar i$. Write $F_n^{\to}$ for the right-moving functions of \cite[equation~(13)]{Chatterjee-Schutz-2010} or \cite[Section~3]{Lee-2020}. Inverting the contour variable in \eqref{eq:Fdef} gives $F_n^{\to}(-s;t)=F_{-n}(s;t)$. Let $d_{ab}$ count adjacent descents between slots $a$ and $b$ in the reflected word. These are the $AB$ pairs of \cite{Chatterjee-Schutz-2010} with $A=2$, $B=1$, hence the reversals of the original $12$ pairs. The index in \cite[Section~III, equation~(28)]{Chatterjee-Schutz-2010} and \cite[Theorem~3.1]{Lee-2020} is therefore
\[
 n(\bar i,\bar j)=\sgn(\bar i-\bar j)\bigl(|i-j|-d_{\bar i\bar j}\bigr)=e_j-e_i .
\]
Reversing both rows and columns gives \eqref{eq:anyword} entry by entry without changing the determinant.
\end{remark}

\subsection{Total crossing}
\begin{theorem}\label{thm:totalcross}
Let $N\ge2$, $p\in[0,1]$, $t\ge0$, and let $X,Y$ be strictly ordered integer configurations of size $N$. With $h(\xi)=1-1/\xi$ and $F_n(x)=F_n(x;t)$,
\begin{align}
P_{(Y,\wa_1)}(X,\wa_N;t)&=\oint_\cC\xi^{\,x_N-y_1-1}e^{\eps(\xi)t}\,\det\Big[F_{i-j}(x_i-y_{j+1})-h(\xi)\,F_{i-j-1}(x_i-y_{j+1})\Big]_{i,j=1}^{N-1}d\xi,\label{eq:totA}\\
P_{(Y,\wb_N)}(X,\wb_1;t)&=\oint_\cC\xi^{\,x_1-y_N-1}e^{\eps(\xi)t}\,\det\Big[F_{i-j}(x_{i+1}-y_{j})-\frac{F_{i-j+1}(x_{i+1}-y_{j})}{h(\xi)}\Big]_{i,j=1}^{N-1}d\xi.\label{eq:totB}
\end{align}
\end{theorem}

\begin{proof}
We prove \eqref{eq:totB}. By Remark~\ref{rem:special}(ii) only permutations with $\sigma(1)=N$ contribute. Write $\varrho(i):=\sigma(i+1)$, $1\le i\le N-1$, so that $\varrho=\sigma\setm N\in S_{N-1}$. By \eqref{eq:factor} and Lemma~\ref{lem:Sprod},
\[
[\bA_\sigma]_{\wb_1,\wb_N}=\sgn(\varrho)\prod_{i=1}^{N-1}h_{\varrho(i)}^{\,i-\varrho(i)}\Big(1-\frac{h_{\varrho(i)}}{h_N}\Big),
\]
and the monomial in \eqref{eq:Lee} is $\xi_N^{\,x_1-y_N-1}\prod_{i=1}^{N-1}\xi_{\varrho(i)}^{\,x_{i+1}-y_{\varrho(i)}-1}$. For fixed $\xi_N$, integrating over $\xi_1,\dots,\xi_{N-1}$ as in the proof of Theorem~\ref{thm:nocross} gives the determinant in \eqref{eq:totB} with $\xi=\xi_N$. The proof of \eqref{eq:totA} is identical: only permutations with $\sigma(N)=1$ contribute, $\varrho:=\sigma\setm1$ is a word in $2,\dots,N$ with $\operatorname{rk}_\varrho(v)=v-1$, and $[\bA_\sigma]_{\wa_N,\wa_1}=\sgn(\varrho)\prod_{i=1}^{N-1}h_{\varrho(i)}^{\,i-\varrho(i)+1}(1-h_1/h_{\varrho(i)})$; one then sets $j=\varrho(i)-1$.
\end{proof}

\begin{remark}[Comparison with the TASEP formula]\label{rem:dgmw}
At $p=0$, \eqref{eq:totA} and \eqref{eq:totB} are, after reflection, the cases $m=N-1$ and $m=1$ of \cite[Proposition~2.1, equation~(2.11)]{deGier-Mead-Wheeler-2023}, where $m$ is the number of particles of type~$2$. In their notation take $\mu_i=-y_{N+1-i}$ and $\nu_i=-x_{N+1-i}$. Each cross factor is $(w-z)/(1-z)=1-(1-w)/(1-z)$. For $m=N-1$ put $z_i=\xi_{N+1-i}^{-1}$ and $w_1=\xi_1^{-1}$; the factor is $1-h_1/h_{N+1-i}$. For $m=1$ put $z_1=\xi_N^{-1}$ and $w_i=\xi_{N-i}^{-1}$; it is $1-h_{N-i}/h_N$. Reverse both indices of the nonsingleton determinant and transpose it, so its columns correspond to the initial sites. Contour inversion gives the measure $\xi^{\,x_i-y_j-1}e^{t(\xi-1)}d\xi$ for the corresponding pair of sites. Integrating each column with its cross factor gives the two differences of $F$'s in \eqref{eq:totA} and \eqref{eq:totB}, respectively. For $p>0$ only the exponential factor changes.
\end{remark}

\begin{remark}[Other final orders]\label{rem:general}
For $1<m<N$, group the permutations by the slot $r$ of the distinguished spectral letter. The admissible slots are $m\le r\le N$ for $\wa_m$ by \eqref{eq:matA}, and $1\le r\le m$ for $\wb_m$ by \eqref{eq:matB}: there are $N-m+1$ and $m$ groups, respectively. For fixed $r$ and fixed distinguished variable, each extra $S$ or $T$ factor depends on just one of the remaining variables and its slot. Lemma~\ref{lem:Sprod} supplies the sign and separable powers for the remaining permutation. Its signed sum is thus an $(N-1)\times(N-1)$ determinant, leaving one outer integral per group. At $p=0$ these probabilities are contained in the general formula of \cite[Theorem~2.1]{deGier-Mead-Wheeler-2023}. We do not use them below and omit the formulas.
\end{remark}

\section{Further results: crossing probabilities and projections}\label{sec:marg}
The $p=1$ no-crossing reduction from $\wa_1$ and total-crossing reduction from $\wb_N$ originate in the 2026 MSc work of Alikhan Aliaidar and Tamirlan Mukazhanov, respectively, at Nazarbayev University. We derive the normalized formulas below for all $p\in[0,1]$. In the two marginals $\FF_N$ and $\GG_N$, the initial sites enter only through $(y_1,y_2)$ in \eqref{eq:FN} and $(y_1,y_N)$ in \eqref{eq:GN}, respectively. Theorems~\ref{thm:redA} and~\ref{thm:redB} give the corresponding reductions for $2\le k\le N$.

Let $\wa_m^{(k)},\wb_m^{(k)}$ be the length-$k$ versions of $\wa_m,\wb_m$, and let $P^{(k)}$ denote the $k$-particle kernel. The following reductions are identities for joint masses of endpoints and crossing events, not conditional endpoint laws.

\subsection{Initial order \texorpdfstring{$12\cdots2$}{12...2}: no crossing}
\begin{theorem}\label{thm:redA}
Let $2\le k\le N$ and $x_1<\dots<x_k$. Then
\begin{equation}\label{eq:redA}
\sum_{x_k<x_{k+1}<\dots<x_N}P_{(Y,\wa_1)}\big((x_1,\dots,x_N),\wa_1;t\big)=P^{(k)}_{((y_1,\dots,y_k),\wa^{(k)}_1)}\big((x_1,\dots,x_k),\wa^{(k)}_1;t\big).
\end{equation}
In particular, the left-hand side does not depend on $y_{k+1},\dots,y_N$. Consequently, for every $N\ge2$ and $x\in\Z$,
\begin{equation}\label{eq:FN}
\begin{aligned}
\FF_N(x;t)&:=\sum_{x<x_2<\dots<x_N}P_{(Y,\wa_1)}\big((x,x_2,\dots,x_N),\wa_1;t\big)\\
&\;=F_0(x-y_1;t)\big[1-F_{-1}(x-y_2;t)\big]-F_0(x-y_2;t)\big[1-F_{-1}(x-y_1;t)\big].
\end{aligned}
\end{equation}
\end{theorem}

\begin{proof}
It suffices to prove \eqref{eq:redA} for $k=N-1\ge2$ and to iterate. In the determinant \eqref{eq:detA} only the last row depends on $x_N$, and its entries are $F_{N-e^{\wa}_j}(x_N-y_j)$, where $N-e^{\wa}_j\ge1$ for $j<N$ and $N-e^{\wa}_N=0$. By Lemma~\ref{lem:F}(d),(e) the sums over $x_N>x_{N-1}$ of these entries converge absolutely and are equal to
\[
-F_{N-1-e^{\wa}_j}(x_{N-1}-y_j)\quad(j<N),\qquad 1-F_{-1}(x_{N-1}-y_N)\quad(j=N).
\]
Since $N\ge3$ we have $e^{\wa}_{N-1}=N-1$, so the row $N-1$ of the matrix is $\big(F_{N-1-e^{\wa}_j}(x_{N-1}-y_j)\big)_{j=1}^N$, whose last entry is $F_{-1}(x_{N-1}-y_N)$. Hence the summed last row equals $(0,\dots,0,1)$ minus the row $N-1$. By multilinearity, the sum over $x_N$ of the determinant is the leading principal minor of order $N-1$, which is the determinant \eqref{eq:detA} for the $(N-1)$-particle system with initial positions $(y_1,\dots,y_{N-1})$. Finally, for $N=2$, summing the second row of $\det\big[F_0(x_i-y_j)\big]_{i,j=1}^2$ over $x_2>x_1=x$ by Lemma~\ref{lem:F}(e) gives \eqref{eq:FN}.
\end{proof}

\subsection{Initial order \texorpdfstring{$1\cdots12$}{1...12}: total crossing}
\begin{theorem}\label{thm:redB}
Let $2\le k\le N$ and $x_1<\dots<x_k$. Then
\begin{equation}\label{eq:redB}
\sum_{x_k<x_{k+1}<\dots<x_N}P_{(Y,\wb_N)}\big((x_1,\dots,x_N),\wb_1;t\big)=P^{(k)}_{((y_1,\dots,y_{k-1},y_N),\wb^{(k)}_k)}\big((x_1,\dots,x_k),\wb^{(k)}_1;t\big).
\end{equation}
In particular, the left-hand side does not depend on $y_k,\dots,y_{N-1}$. Consequently, for every $N\ge2$ and $z<x$,
\begin{equation}\label{eq:GNjoint}
\begin{aligned}
\sum_{x<x_3<\dots<x_N}&P_{(Y,\wb_N)}\big((z,x,x_3,\dots,x_N),\wb_1;t\big)\\
&=F_0(z-y_N;t)\,F_0(x-y_1;t)-F_{-1}(z-y_N;t)\,F_1(x-y_1;t),
\end{aligned}
\end{equation}
and
\begin{equation}\label{eq:GN}
\begin{aligned}
\GG_N(x;t)&:=\sum_{z<x}\ \sum_{x<x_3<\dots<x_N}P_{(Y,\wb_N)}\big((z,x,x_3,\dots,x_N),\wb_1;t\big)\\
&\;=F_0(x-y_1;t)\,F_{-1}(x-1-y_N;t)-F_1(x-y_1;t)\,F_{-2}(x-1-y_N;t).
\end{aligned}
\end{equation}
\end{theorem}

\begin{proof}
For $N\ge3$ write $g_{ij}=F_{i-j}(x_{i+1}-y_j)-h(\xi)^{-1}F_{i-j+1}(x_{i+1}-y_j)$ in \eqref{eq:totB}. By Lemma~\ref{lem:F}(e),
\[
 \sum_{x_N>x_{N-1}}g_{N-1,j}=\mathbf1_{\{j=N-1\}}-g_{N-2,j}.
\]
The sum is absolutely convergent uniformly on $\cC$, since $|h(\xi)^{-1}|\le(1-1/R)^{-1}$ and the last-row indices are nonnegative. Taking it inside the integral leaves the leading minor and the same outer factor involving $y_N$. This removes $y_{N-1}$; iteration gives \eqref{eq:redB}. For $k=2$ the remaining integral gives \eqref{eq:GNjoint}, since the outer factors $1$ and $h^{-1}$ integrate to $F_0(z-y_N)$ and $F_{-1}(z-y_N)$. Summation over $z\le x-1$ gives \eqref{eq:GN}.
\end{proof}

\begin{remark}[Double-integral forms and pairwise cancellation]\label{rem:theses}
By \eqref{eq:factor}, \eqref{eq:FN} and \eqref{eq:GN} are equivalent to
\begin{align*}
\FF_N(x;t)&=\oint_\cC\!\oint_\cC T_{21}\,\frac{\xi_1}{\xi_1-1}\prod_{k=1}^{2}\xi_k^{\,x-y_k-1}e^{\eps(\xi_k)t}\,d\xi_1d\xi_2+F_0(x-y_1;t)-F_0(x-y_2;t),\\
\GG_N(x;t)&=\oint_\cC\!\oint_\cC \frac{T_{21}}{\xi_2-1}\ \xi_1^{\,x-y_1-1}\xi_2^{\,x-y_N-1}e^{\eps(\xi_1)t+\eps(\xi_2)t}\,d\xi_1d\xi_2,
\end{align*}
respectively. Both MSc arguments at $p=1$ pair permutations differing in their last two letters. In the normalization used here, summing over the largest gap gives the cancellation $S_{\beta\alpha}h_\alpha^{-1}=-h_\beta^{-1}$.
\end{remark}

\subsection{Crossing probabilities}
Let $I_k$ denote the modified Bessel function of the first kind, $e^{t(z+1/z)}=\sum_{k\in\Z}I_k(2t)z^k$.

\begin{corollary}\label{cor:cross}
Let $N\ge2$, and write $d_{12}=y_2-y_1$ and $d_{1N}=y_N-y_1$.
\begin{enumerate}[label=\textup{(\alph*)},leftmargin=2.2em]
\item Started from $(Y,\wa_1)$, the probability that the order of the species at time $t$ is still $12\cdots2$ is
\[
\sum_{x\in\Z}\FF_N(x;t)=e^{-2t}\sum_{k=1-d_{12}}^{d_{12}}I_k(2t).
\]
\item Started from $(Y,\wb_N)$, the probability that the order of the species at time $t$ is $21\cdots1$ is
\[
\sum_{x\in\Z}\GG_N(x;t)=1-e^{-2t}\sum_{k=1-d_{1N}}^{d_{1N}}I_k(2t).
\]
\end{enumerate}
Both probabilities are independent of $p$. For each fixed integer $d\ge1$, $e^{-2t}\sum_{k=1-d}^{d}I_k(2t)= d\,(\pi t)^{-1/2}\,(1+o(1))$ as $t\to\infty$.
\end{corollary}

\begin{proof}
(a) Let $W,W'$ be independent walks started at zero. By \eqref{eq:FN} and Lemma~\ref{lem:F}(c),
\[
\sum_x\FF_N(x;t)=\PP(W_t+y_1<W'_t+y_2)-\PP(W'_t+y_2<W_t+y_1)=\PP(-d_{12}<\Delta_t\le d_{12}),
\]
where $\Delta_t:=W'_t-W_t$ is symmetric. Since $\Delta$ jumps by $\pm1$ at rate $p+q=1$ each, $\EE z^{\Delta_t}=e^{-2t}e^{t(z+1/z)}$ and $\PP(\Delta_t=k)=e^{-2t}I_k(2t)$, whatever the value of $p$. (b) By \eqref{eq:redB} with $k=2$, the left-hand side is the probability that the two-particle system started from $((y_1,y_N),12)$ has order $21$ at time $t$, which is one minus the probability in (a) with initial separation $d_{1N}$. The asymptotic formula follows from the local limit theorem for $\Delta_t$.
\end{proof}

The first exchange for $\wa_1$ and the last exchange for $\wb_N$ each have the law of a two-walk meeting time, irrespective of the other positions or of $p$. In the notation of Section~\ref{subsec:overtaking}, Corollary~\ref{cor:cross}(b) gives $\PP(R(t)=0)$; taking its complement gives $\PP(R(t)\ge1)$, as in Theorem~\ref{thm:overtaking} for $m=1$.

\subsection{Probabilistic explanation}\label{sec:prob}
The reductions of Theorems~\ref{thm:redA} and~\ref{thm:redB} reflect exact projection properties of the dynamics, which we now describe. They give an independent derivation of the fact that the left-hand sides of \eqref{eq:redA} and \eqref{eq:redB} are $k$-particle quantities, but not of the formulas of Section~\ref{sec:trans}.

\begin{proposition}\label{prop:projA}
Start from $(Y,\wa_1)$. Here $R(0)=1$: the first exchange occurs when the sole second-class particle ceases to be leftmost, at $\tau_1$ of Section~\ref{subsec:overtaking}. Retain the second-class particle and the $k-1$ leftmost first-class particles, $2\le k\le N$, and send their positions to a cemetery state at $\tau_1$. This projected chain has the law of the $k$-particle system from $((y_1,\ldots,y_k),\wa^{(k)}_1)$ killed at its own first exchange. For $k=2$ it is the law of two independent walks from $y_1<y_2$ killed at their first meeting.
\end{proposition}

\begin{proof}
Before $\tau_1$ the second-class particle is leftmost. Its left jumps are unrestricted. A right attempt leads to the adjacent site if vacant and causes an exchange otherwise. The first-class particles form a single-species PushASEP (Remark~\ref{rem:proj}): a right attempt pushes only to the right, and a left attempt can be blocked only by the left neighbour. A left attempt by the leftmost first-class particle causes an exchange if the second-class particle is its left neighbour. Thus the surviving transition rates and the killing rate depend only on the retained positions and agree with those of the $k$-particle system. For $k=2$ the two walks are independent until one attempts to jump onto the other, at which point the projection is killed.
\end{proof}

The word at time $t$ is $\wa_1$ exactly on $\{\tau_1>t\}$. Thus Proposition~\ref{prop:projA} gives \eqref{eq:redA}; for $k=2$ the Karlin--McGregor formula \cite{Karlin-McGregor-1959} gives \eqref{eq:FN}.

\begin{proposition}\label{prop:projB}
Consider a system with a single first-class particle and at least $k-1$ second-class particles, $k\ge2$, in an arbitrary initial state. Let $Z(t)$ be the position of the first-class particle and $L_1(t)<\dots<L_{k-1}(t)$ the positions of the $k-1$ leftmost second-class particles at time $t$. Then $(L_1,\dots,L_{k-1};Z)$ is a Markov chain with the law of the $k$-particle system with one first-class particle.
\end{proposition}

\begin{proof}
Let $\eta:=\{L_1,\dots,L_{k-1},Z\}$ be the retained set of occupied sites, with $Z$ marked as first-class. We compare the changes of the sites and of the mark with those in the $k$-particle system. Call the remaining second-class particles \emph{outer}; they all lie to the right of $L_{k-1}$.

\emph{Right jump of the particle at $L_i$.} In the $k$-particle system it lands on the nearest site $b>L_i$ not in $\eta$. All sites strictly between $L_i$ and $b$ belong to $\eta$. If $b$ is empty in the full system, the particle lands on $b$ there as well, and $b$ lies to the left of all outer particles. If $b$ is occupied in the full system, its occupant is the leftmost outer particle; the jumping particle lands further to the right, and in the set of the $k-1$ leftmost second-class particles the site $L_i$ is replaced by $b$. In both cases $\eta$ changes as in the $k$-particle system.

\emph{Right jump of the first-class particle.} It moves to $Z+1$. If $Z+1=L_i$, the displaced particle moves to the nearest empty site and the previous argument applies. If $Z+1$ is occupied by an outer particle, this particle is displaced to the right and $\eta$ changes as in the $k$-particle system, where $Z+1$ is empty.

\emph{Left jumps.} The particle at $L_i$ moves to $L_i-1$ unless this site is occupied; an occupant can only be the particle at $L_{i-1}$ or the first-class particle, both in $\eta$. The first-class particle always moves to $Z-1$, exchanging with a second-class particle if there is one; if the latter is outer, it moves to $Z$, which is to the right of $L_{k-1}$, and $L_1,\dots,L_{k-1}$ are unchanged.

\emph{Jumps of outer particles} do not change $\eta$: to the right they move away, and to the left they are suppressed by, or stay to the right of, the particle at $L_{k-1}$ or the first-class particle.
\end{proof}

For $(Y,\wb_N)$ this projection is the $k$-particle system from $((y_1,\ldots,y_{k-1},y_N),\wb^{(k)}_k)$. Since the full order is $\wb_1$ exactly when $Z(t)<L_1(t)$, this proves \eqref{eq:redB} and gives the analogous projection for every relative rank of $Z$ among the retained particles. For $k=2$, while $L_1<Z$ the killed description agrees with two independent walks killed at their first meeting. The physical chain exchanges at that jump; afterwards $Z<L_1$ and the first-class particle can push the second-class particle.

\begin{remark}[Rank selection and recolouring]\label{rem:recolour}
Proposition~\ref{prop:projB} follows particles by their current rank among the second-class particles, not by permanent tags. It is related to the monotone-recolouring property of multi-species push systems \cite[Proposition~2.2]{AMW}, where colours are merged or erased by a weakly increasing map fixing $0$.
\end{remark}

\section{Proofs of the persistence results}\label{sec:persist-proofs}

\subsection{Proof of Theorem~\ref{thm:sN}}\label{subsec:sN-proof}
Summing a Sch\"utz-type transition determinant over ordered particle positions is a classical route to current distributions; cf.\ \cite[Section~II]{Rakos-Schutz-2005}. We use the following telescoping identity. For a function $f$ on $\Z$ let $(\nabla f)(x):=f(x)-f(x-1)$. For row vectors $v_1,\dots,v_n$ of length $n$ we write $\det(v_1,\dots,v_n)$ for the determinant of the matrix with rows $v_1,\dots,v_n$.

\begin{lemma}\label{lem:sum}
Let $R>1$ and $m\ge n\ge1$, and let $G_1,\dots,G_m$ be complex functions on $\Z$ such that $R^{-x}|G_j(x)|$ is bounded for $x\le0$, for every $j$. For a set $J\subset\{1,\dots,m\}$ of $n$ elements let $g_J(x):=(G_j(x))_{j\in J}$. Then, for every $b\in\Z$,
\[
\sum_{x_1<\dots<x_n\le b}\det\big(\nabla g_J(x_1),\nabla^2g_J(x_2),\dots,\nabla^ng_J(x_n)\big)=\det\big(g_J(b-n+1),\dots,g_J(b)\big),
\]
the series converging absolutely.
\end{lemma}

The proof, by induction on $n$, is given in Appendix~\ref{app:lemmas}.

\begin{proof}[Proof of Theorem~\ref{thm:sN}]
For $1\le j\le N$ let
\[
G_j(x):=F_{-e^{\wb}_j}(x-y_j;t),\qquad g(x):=\big(G_1(x),\dots,G_N(x)\big).
\]
By Lemma~\ref{lem:F}(d) the hypothesis of Lemma~\ref{lem:sum} holds, and by Lemma~\ref{lem:F}(b) the function $\nabla^nG_j$ is $x\mapsto F_{n-e^{\wb}_j}(x-y_j;t)$, for every $n\ge0$. Thus the rows of the matrix in \eqref{eq:detB} are $\nabla^ig(x_i)$ for $i\le N-1$, and $\nabla^{N-1}g(x_N)$ for $i=N$. Expanding along the last row and applying Lemma~\ref{lem:sum} with $n=N-1$ and $b=x-1$ to the cofactors, we get
\begin{align*}
s_N(x;t)&=\det\big(g(x-N+1),\dots,g(x-1),\nabla^{N-1}g(x)\big)\\
&=\det\big(g(x-N+1),\dots,g(x-1),g(x)\big),
\end{align*}
because $\nabla^{N-1}g(x)-g(x)$ is a linear combination of the other rows. Finally, $\nabla^kg(x)$ is $(-1)^kg(x-k)$ plus a linear combination of $g(x),\dots,g(x-k+1)$, so that
\begin{align*}
\det\big(g(x),\nabla g(x),\dots,\nabla^{N-1}g(x)\big)&=(-1)^{N(N-1)/2}\det\big(g(x),g(x-1),\dots,g(x-N+1)\big)\\
&=\det\big(g(x-N+1),\dots,g(x)\big),
\end{align*}
and $(\nabla^{i-1}G_j)(x)=F_{i-1-e^{\wb}_j}(x-y_j;t)$.
\end{proof}

\subsection{Proofs of Theorem~\ref{thm:decay} and Corollaries~\ref{cor:harmonic} and~\ref{cor:doob}}\label{subsec:decay-proof}
We first remove the near-degeneracy of the determinant \eqref{eq:sN}, whose last two columns differ only through $y_{N-1}$ and $y_N$.

\begin{lemma}\label{lem:colred}
For $k\in\Z$ let
\[
D(x;k):=\det\Big[\,F_{i-1-j}(x-y_j;t)\ \ (1\le j\le N-1)\ \ \Big|\ \ F_{i-1}(x-k;t)\,\Big]_{i=1}^N
\]
be the determinant whose first $N-1$ columns are those of \eqref{eq:sN} and whose last column is $(F_{i-1}(x-k;t))_{i=1}^N$. Then
\[
s_N(x;t)=(-1)^{N-1}\sum_{(k_1,\dots,k_{N-1})}D(x;k_{N-1}),
\]
where the sum is over the $\gamma_N(Y)$ sequences satisfying \eqref{eq:chains}.
\end{lemma}

\begin{proof}
Let $c^{(l)}(k)$ denote the column $(F_{i-1-l}(x-k;t))_{i=1}^N$. By Lemma~\ref{lem:F}(b), for integers $a<b$,
\[
c^{(l)}(b)-c^{(l)}(a)=-\sum_{k=a}^{b-1}c^{(l-1)}(k).
\]
The columns of \eqref{eq:sN} are $c^{(1)}(y_1),\dots,c^{(N-1)}(y_{N-1})$ and $c^{(N-1)}(y_N)$. Subtracting the column $N-1$ from the last column replaces the latter by $-\sum_{k_1=y_{N-1}}^{y_N-1}c^{(N-2)}(k_1)$. In each term we subtract the column $N-2$, which replaces $c^{(N-2)}(k_1)$ by $-\sum_{k_2=y_{N-2}}^{k_1-1}c^{(N-3)}(k_2)$, and so on. After $N-1$ steps the last column is $c^{(0)}(k_{N-1})$ and the accumulated sign is $(-1)^{N-1}$.
\end{proof}

Next we need estimates for the functions $F_n$ on the diffusive scale: a local limit theorem together with bounds which are uniform in $x$. For $u\in\mathbb R$ put $u^\pm:=\max(\pm u,0)$. The functions $\beta_n(\cdot)$ below are envelopes, distinct from the exponents in Section~\ref{subsec:capture}:
\[
\beta_n(u):=e^{-|u|}\ \ (n\ge0),\qquad \beta_n(u):=e^{-u^-}(1+u^+)^{-n-1}\ \ (n\le-1).
\]

\begin{lemma}\label{lem:llt}
Let $n\in\Z$. There is a constant $C_n$ depending only on $n$ such that, for all $p\in[0,1]$, $t\ge1$ and $x\in\Z$, with $u=(x-vt)/\sqrt t$,
\begin{equation}\label{eq:lltbound}
\big|t^{(n+1)/2}F_n(x;t)\big|\le C_n\,\beta_n(u),
\end{equation}
and for every $A>0$,
\begin{equation}\label{eq:lltlimit}
\lim_{t\to\infty}\ \sup_{x\in\Z:\ |x-vt|\le A\sqrt t}\Big|t^{(n+1)/2}F_n(x;t)-\Phi_{-n}\Big(\frac{x-vt}{\sqrt t}\Big)\Big|=0 .
\end{equation}
\end{lemma}

The proof is given in Appendix~\ref{app:lemmas}.

\begin{proof}[Proof of Theorem~\ref{thm:decay}]
(a) For $l=0$ the identity below follows by integration from $-\infty$. For $l>0$, differentiate in $u$ and apply the identity at $l-1$; both sides vanish as $u\to-\infty$ by Gaussian decay. Thus for all $l,m\ge0$,
\[
\int_{-\infty}^u\varphi^{(m)}(r)\,(u-r)^l\,dr=l!\,\Phi_{l+1-m}(u).
\]
Hence the entries of the first column of $\omega_N$ are $\Phi_{1-i}(u)=\varphi^{(i-1)}(u)$, and those of the column $j\ge2$ are
\[
\Phi_{j-i}(u)=\frac1{(j-2)!}\int_{-\infty}^u\varphi^{(i-1)}(r)\,(u-r)^{j-2}\,dr .
\]
By multilinearity in the columns $2,\dots,N$, with $z:=(u,r_2,\dots,r_N)$,
\[
\omega_N(u)=\int_{(-\infty,u]^{N-1}}\det\big[\varphi^{(i-1)}(z_k)\big]_{i,k=1}^N\ \prod_{j=2}^N\frac{(u-r_j)^{j-2}}{(j-2)!}\ dr_2\cdots dr_N .
\]
Since $\varphi^{(i-1)}=(-1)^{i-1}\mathrm{He}_{i-1}\varphi$ with a monic polynomial $\mathrm{He}_{i-1}$ of degree $i-1$, the determinant equals $(-1)^{N(N-1)/2}\Delta_N(z)\prod_k\varphi(z_k)$; it is antisymmetric in $(r_2,\dots,r_N)$. We may therefore replace the product over $j$ by its antisymmetrization,
\[
\frac1{(N-1)!}\det\Big[\frac{(u-r_k)^{j-2}}{(j-2)!}\Big]_{j,k=2}^N=\frac{(-1)^{(N-1)(N-2)/2}}{(N-1)!\,\prod_{l=1}^{N-2}l!}\,\Delta_{N-1}(r_2,\dots,r_N).
\]
Together with $\Delta_N(z)=(-1)^{N-1}\prod_j(u-r_j)\,\Delta_{N-1}(r_2,\dots,r_N)$, and since the three signs multiply to $(-1)^{N(N-1)}=1$, this gives
\[
\omega_N(u)=\frac{\varphi(u)}{(N-1)!\,\prod_{l=1}^{N-2}l!}\int_{(-\infty,u]^{N-1}}\Delta_{N-1}(r)^2\prod_{j}(u-r_j)\prod_j\varphi(r_j)\,dr ,
\]
which is \eqref{eq:Wformula} because the integrand is symmetric and $\Delta_{N-1}(r)\prod_j(u-r_j)=\Delta_N(r,u)$. Integrating over $u$ gives \eqref{eq:kappa}; the integral is finite and the integrand is positive on the domain of integration.

(b) Fix $k\in\Z$ and expand
\[
D(x;k)=\sum_{\sigma\in S_N}\sgn(\sigma)\prod_{j=1}^{N-1}F_{\sigma(j)-1-j}(x-y_j;t)\cdot F_{\sigma(N)-1}(x-k;t).
\]
In each term the indices $n_j:=\sigma(j)-1-j$ $(j<N)$ and $n_N:=\sigma(N)-1\ge0$ satisfy $\sum_{j=1}^N(n_j+1)=N$. Put $(s_1,\ldots,s_N)=(y_1,\ldots,y_{N-1},k)$, $s_*=\min_j s_j$, $d_*=\max_j s_j-s_*$, and $u_*=(x-s_*-vt)/\sqrt t$. For $a\ge0$, $\beta_n(w-a)\le e^a\beta_n(w)$: for $n<0$ the envelope is nondecreasing, and for $n\ge0$ use the triangle inequality. Applying \eqref{eq:lltbound} at $u_*-(s_j-s_*)/\sqrt t$ in column $j$ gives
\begin{equation}\label{eq:Dbound}
 \big|t^{N/2}D(x;k)\big|\le C_N e^{N d_*/\sqrt t}(1+|u_*|)^{N^2}e^{-|u_*|}\qquad(t\ge1,\ x\in\Z).
\end{equation}
Here $C_N$ depends only on $N$: the last column supplies the exponential and the other columns contribute at most the displayed polynomial. For fixed shifts, with $u=(x-vt)/\sqrt t$, \eqref{eq:lltlimit} gives, uniformly for $|u|\le A$,
\[
t^{N/2}D(x;k)\ \longrightarrow\ \det\Big[\,\Phi_{j-i+1}(u)\ \ (1\le j\le N-1)\ \ \Big|\ \ \Phi_{1-i}(u)\,\Big]_{i=1}^N=(-1)^{N-1}\omega_N(u),
\]
the last equality being obtained by moving the last column to the first place. With Lemma~\ref{lem:colred} this proves (b).

(c), (d) For an absolutely summable function $f$ on $\Z$ we have $\sum_{x\in\Z}f(x)=\sqrt t\int_{\mathbb R}f(\lceil vt+r\sqrt t\,\rceil)\,dr$. Hence, for $a<b$,
\[
t^{(N-1)/2}\sum_{x:\ a<(x-vt)/\sqrt t\le b}s_N(x;t)=\int_{a_t}^{b_t}t^{N/2}s_N\big(\lceil vt+r\sqrt t\,\rceil;t\big)\,dr\ \longrightarrow\ \gamma_N(Y)\int_a^b\omega_N(u)\,du
\]
by (b), where $a_t\to a$ and $b_t\to b$. By \eqref{eq:Dbound} and Lemma~\ref{lem:colred} the integrand is bounded by $C(1+|r|)^{N^2}e^{-|r|}$ uniformly in $t\ge1$, so that the same holds with $a=-\infty$ and $b=\infty$, which is (c). Dividing the two limits gives (d).
\end{proof}

\begin{proof}[Proof of Corollary~\ref{cor:harmonic}]
By the Markov property at time $s$,
\[
\mathcal S_N(t+s;Y)=\EE_Y\big[\mathcal S_N(t;Y(s));\ \tau>s\big].
\]
Multiply by $t^{(N-1)/2}$. For each surviving state $Y(s)$ the integrand converges pointwise by \eqref{eq:decay}. For an ordered initial state $Y'$, translate $y'_1$ to zero and let $d=y'_N-y'_1$. All column shifts in Lemma~\ref{lem:colred} then lie in $[0,d]$, with $s_*=0$, so the shift factor in \eqref{eq:Dbound} is at most $e^{Nd}$. The sum of the polynomial--exponential envelope over $x$, times $t^{-1/2}$, is bounded independently of $t\ge1$ and of the lattice translate: on each cell of length $t^{-1/2}\le1$ its endpoint value is bounded by a constant times its value anywhere in the cell. Summing over the $\gamma_N(Y')\le d^{N-1}$ chains gives
\begin{equation}\label{eq:survival-uniform-Y}
t^{(N-1)/2}\,\mathcal S_N(t;Y')\le C'\,\gamma_N(Y')\,e^{Nd}\le C''e^{(N+1)d},\qquad t\ge1.
\end{equation}
On $\{\tau>s\}$ the diameter is at most $y_N-y_1+J_s$: a surviving update changes either extreme rank by at most one and cannot move both extremes. Here $J_s$ is Poisson with mean $Ns$ and has every exponential moment. The bound \eqref{eq:survival-uniform-Y} permits dominated convergence; cancel $\kappa_N>0$.
\end{proof}

\begin{proof}[Proof of Corollary~\ref{cor:doob}]
Since $\gamma_N>0$, the kernel \eqref{eq:doob-kernel} is well defined. It inherits the semigroup property from $K_t$, and $\sum_{X}K_t^\gamma(Y,X)=\EE_Y[\gamma_N(Y(t));\,\tau>t]/\gamma_N(Y)=1$ by Corollary~\ref{cor:harmonic}. By the same corollary, $\mathbf1_{\{\tau>s\}}\gamma_N(Y(s))/\gamma_N(Y)$ is a nonnegative martingale with mean one, so that the right-hand side of \eqref{eq:long-conditioning} is the expectation of $\Psi$ for the process with kernel $K^\gamma$. For $\mathcal T>s$, the Markov property at time $s$ gives
\[
\EE_Y[\Psi\mid\tau>\mathcal T]
=\EE_Y\left[\Psi\,\mathbf1_{\{\tau>s\}}\,
       \frac{\mathcal S_N(\mathcal T-s;Y(s))}{\mathcal S_N(\mathcal T;Y)}\right].
\]
On $\{\tau>s\}$, Theorem~\ref{thm:decay}(c) gives the convergence of the ratio to $\gamma_N(Y(s))/\gamma_N(Y)$ as $\mathcal T\to\infty$. Moreover, $\mathcal T^{(N-1)/2}\mathcal S_N(\mathcal T;Y)$ is bounded away from zero for all sufficiently large $\mathcal T$, and if $\mathcal T\ge2s+2$ then $\mathcal T/(\mathcal T-s)\le2$, so that \eqref{eq:survival-uniform-Y} bounds the ratio by $C_Y\exp\big((N+1)(W(s)-x_1(s))\big)$. This bound is integrable, as in the proof of Corollary~\ref{cor:harmonic}, and dominated convergence proves \eqref{eq:long-conditioning}.
\end{proof}

\subsection{Proof of Proposition~\ref{prop:packed-km}}\label{subsec:packed-proof}
\begin{proof}[Proof of Proposition~\ref{prop:packed-km}]
Put $c_M=\prod_{r=0}^{M-1}r!$, so that $\Delta_M(Y_-)=c_M$, and set
\[
B_i(z)=\binom{x-z-1}{M-i},\qquad 1\le i\le M,\quad z<x.
\]
Their degrees decrease from $M-1$ to $0$, and their leading coefficients give
\[
\det[B_i(z_k)]_{i,k=1}^M=\frac{\Delta_M(z_1,\dots,z_M)}{c_M}.
\]
Expanding along the fixed last row and applying the discrete Cauchy--Binet identity to each minor turns the sum in \eqref{eq:packed-km} into
\[
\det\big[F_{i-N}(x-y_j+i-N;t)\big]_{i,j=1}^N.
\]
Indeed, for $i<N$ the moment in row $i$, column $j$ is
\[
\sum_{z<x}\binom{x-z-1}{M-i}w_t(z-y_j)
       =F_{i-N}(x-y_j+i-N;t),
\]
by repeated summation of $F_0=w_t$ using Lemma~\ref{lem:F}(e), while the last row is $w_t(x-y_j)$. Absolute convergence follows from the exponential moments of the walk, since the weights are polynomials.

For integer $u$, the identity $F_n(u-1)=F_n(u)-F_{n+1}(u)$ gives
\[
F_{i-N}(x-y_j+i-N;t)
 =\sum_{k=i}^N(-1)^{k-i}\binom{N-i}{k-i}F_{k-N}(x-y_j;t).
\]
This row transformation has determinant $1$, so the preceding determinant equals $\det[F_{i-N}(x-y_j;t)]_{i,j=1}^N$. Denote the columns of this latter matrix by $C_j$. For the first $M$ columns, use the simultaneous triangular transformation
\[
C_j\longmapsto\sum_{k=j}^M(-1)^{k-j}\binom{M-j}{k-j}C_k,
\qquad 1\le j\le M,
\]
leaving the last column unchanged. Since $y_k=y_j+k-j$, the resulting entry in column $j\le M$ is
\[
\sum_{k=j}^M(-1)^{k-j}\binom{M-j}{k-j}
 F_{i-N}(x-y_j-k+j;t)=F_{i-1-j}(x-y_j;t).
\]
The last column is $F_{i-N}(x-y_N;t)$. The column transformation also has determinant $1$, and the resulting matrix is precisely \eqref{eq:sN}, with $e^{\wb}=(1,\dots,M,M)$. This proves \eqref{eq:packed-km} for all $p\in[0,1]$.

For the probabilistic interpretation, the determinant in \eqref{eq:packed-km} is the Karlin--McGregor kernel of $N$ independent walks killed at their first collision \cite{Karlin-McGregor-1959}. The function $\Delta_M$ is harmonic for the lower $M$ walks killed on collision: their generator maps it to an alternating polynomial of degree strictly below the minimal alternating degree $M(M-1)/2$, hence to zero, and it vanishes on collision states. Write $z_-=(z_1,\ldots,z_M)$. In the lower Doob transform a jump from $z_-$ to $\widetilde z_-$ has its rate multiplied by $\Delta_M(\widetilde z_-)/\Delta_M(z_-)$, zero at a lower collision. Harmonicity keeps the total lower jump rate equal to $M$, so this process is conservative. Take its independent product with the upper walk, whose rates are unchanged, and kill the product upon an upper--lower encounter. On each surviving path the lower rate ratios telescope to $\Delta_M(z_-)/\Delta_M(Y_-)$ and the holding rates are unchanged before killing. Thus the weighted determinant is exactly the killed product kernel; summing over $z_-$ gives the stated probability.

Finally, the adjacent-particle count following \eqref{eq:chains} gives $\gamma_N(Y)=\binom{y_N-y_1}{M}$, while
\[
\frac{\Delta_N(Y)}{\Delta_M(Y_-)}
 =\prod_{j=1}^M(y_N-y_1-j+1)=M!\binom{y_N-y_1}{M}.
\]
\end{proof}

For fixed ordered $Y$, not necessarily packed, let $A_N(x;t;Y)$ denote the surviving upper marginal of the auxiliary walks, given by the right-hand side of \eqref{eq:packed-km}. The Cauchy--Binet sum and row operations above do not use packing, so
\[
 A_N(x;t;Y)=\frac{c_M}{\Delta_M(Y_-)}\det[F_{i-N}(x-y_j;t)]_{i,j=1}^N .
\]
Put $a_j=y_j-y_1$ and $v_\ell=(F_{i-N+\ell}(x-y_1;t))_{i=1}^N$. Column $j$ is $\sum_{\ell=0}^{a_j}(-1)^\ell\binom{a_j}{\ell}v_\ell$. In the finite Cauchy--Binet expansion choose $0\le\ell_1<\cdots<\ell_N\le a_N$ and put $E=\sum_j\ell_j-\binom N2\ge0$. For any determinant product the indices $n_i=i-N+\ell_{\sigma(i)}$ satisfy
\[
 \sum_{i=1}^N(n_i+1)=N+E,\qquad n_i\ge 1-N.
\]
The column with $\ell_N\ge M$ supplies an index $n_i\ge0$. With $u=(x-y_1-vt)/\sqrt t$, this factor gives exponential decay in \eqref{eq:lltbound}; the others grow at most polynomially. For fixed $Y$ the index sets are finite, so for all of them, $t\ge1$ and $x\in\Z$,
\[
 t^{(N+E)/2}\big|\det[v_{\ell_1},\ldots,v_{\ell_N}]\big|
 \le C_{N,Y}(1+|u|)^{N^2}e^{-|u|}.
\]
Only $(0,\ldots,M)$ has $E=0$. Its binomial minor is $\det[\binom{a_j}{i-1}]=\Delta_N(Y)/\prod_{k=0}^{M}k!$ by leading coefficients. The factors $(-1)^\ell$ cancel the sign of reversing these columns; by \eqref{eq:lltlimit}, $t^{N/2}\det[v_M,\ldots,v_0]-\omega_N(u)\to0$ uniformly on compact sets. After multiplication by $t^{N/2}$ every other minor is bounded by $t^{-1/2}$ times the displayed envelope. The same envelope controls the tails of the Riemann sums of mesh $t^{-1/2}$. Using $\int\omega_N=\kappa_N$ gives
\[
 \sum_x A_N(x;t;Y)\sim
 \kappa_N\,\frac{\Delta_N(Y)}{M!\,\Delta_M(Y_-)}\,t^{-M/2}.
\]

For the joint limit, let $K_t^{\mathrm{aux}}(Y,Z)$ denote the summand in \eqref{eq:packed-km} at $Z=(z_1,\ldots,z_N)$, $z_N=x$. Its killed-product interpretation above does not require packing. Put $u_i=(z_i-y_1-vt)/\sqrt t$ and let $z_-,u_-$ denote the first $M$ coordinates of $Z,u$, respectively. For fixed $Y$, uniformly on compact subsets of the ordered chamber,
\[
 t^{(2N-1)/2}K_t^{\mathrm{aux}}(Y,Z)
 =\frac{\Delta_N(Y)}{\Delta_M(Y_-)\prod_{k=0}^M k!}\,
 \Delta_M(u_-)\Delta_N(u)\prod_{i=1}^N\varphi(u_i)+o(1).
\]
Indeed, expand $w_t(z_i-y_j)=\sum_{\ell=0}^{a_j}(-1)^\ell\binom{a_j}{\ell}F_\ell(z_i-y_1;t)$ in the determinant. By Lemma~\ref{lem:llt} the minor with indices $\ell_1<\cdots<\ell_N$ is $O(t^{-(N+\sum_j\ell_j)/2})$ on compact sets. Only $0,\ldots,M$ contribute at leading order, with the binomial minor computed above. For each such column, $(-1)^\ell\Phi_{-\ell}=(-1)^\ell\varphi^{(\ell)}$ is a monic polynomial of degree $\ell$ times $\varphi$, giving $\Delta_N(u)\prod_i\varphi(u_i)$. Finally $\Delta_M(z_-)=t^{M(M-1)/4}\Delta_M(u_-)$. By \eqref{eq:kappa} the displayed profile integrates to $\kappa_N\Delta_N(Y)/(M!\,\Delta_M(Y_-))$, the total-mass constant just proved. Compactly supported Riemann sums and equality of the limiting masses therefore imply weak convergence after conditioning on survival; no mass is lost at the boundary or at infinity. The resulting density is the integrand in \eqref{eq:kappa} normalized to unit mass.

\begin{remark}[Initial-coordinate differences]\label{rem:initial-differences}
Let $E_j$ increase the initial coordinate $y_j$ by one, keeping the others fixed, and set $\mathsf D_Y=\prod_{j=1}^{M}(1-E_j)^{M-j}$. Define
\[
 B_N(x;t;Y):=\det[F_{i-N}(x-y_j;t)]_{i,j=1}^N .
\]
For ordered $Y$ this is $\Delta_M(Y_-)A_N(x;t;Y)/c_M$; for shifted configurations we use the determinant, which is alternating in the initial coordinates. Then
\[
 s_N(x;t)=\mathsf D_Y B_N(x;t;Y),\qquad
 \gamma_N(Y)=\frac{\mathsf D_Y\Delta_N(Y)}{M!\,c_M}.
\]
Indeed, $(1-E_j)F_n(x-y_j;t)=F_{n+1}(x-y_j;t)$. In column $j\le M$ the index $i-N$ therefore becomes $i-1-j$; the last column is unchanged. This gives \eqref{eq:sN}. Every shifted lower coordinate is at most $y_j+M-j\le y_M<y_N$. Repeated columns give zero; otherwise reorder the lower coordinates with the permutation sign. For fixed $Y$ sum over $x$ and apply the auxiliary survival asymptotic above at these finitely many starts. Dividing by $\kappa_N t^{-M/2}$ and using Theorem~\ref{thm:decay}(c) gives the second identity. For $Y=(0,2,5)$ it reads $(30-12)/2=9$. For moving initial sites the signed differences can cancel leading terms; the uniform estimate in Theorem~\ref{thm:joint} instead uses the positive chain expansion \eqref{eq:joint-chains}.
\end{remark}

\subsection{Proof of Theorem~\ref{thm:joint}}\label{subsec:joint-proof}
\begin{proof}
For each integer $a$ and column index $\ell$ let $C_\ell(a)$ have $i$th entry $F_{e_i-\ell}(x_i-a;t)$. Equation~\eqref{eq:detB} writes $K_t$ as the determinant with columns $C_1(y_1),\ldots,C_M(y_M),C_M(y_N)$. The telescoping identity
\[
 C_\ell(b)-C_\ell(a)=-\sum_{r=a}^{b-1}C_{\ell-1}(r),\qquad a<b,
\]
holds entrywise even though the endpoint $x_i$ now depends on the row. Subtract the preceding column from the last as in Lemma~\ref{lem:colred}, and repeat $M$ times. Moving the resulting last column to the front cancels the sign $(-1)^M$. Thus
\begin{equation}\label{eq:joint-chains}
 K_t(Y,X)=\sum_{\boldsymbol{k}\ \mathrm{in}\ \eqref{eq:chains}}
 \det\big[F_{e_i-j}(x_i-b_j;t)\big]_{\substack{1\le i\le N\\0\le j\le M}},
 \qquad b_0=k_M,\quad b_j=y_j\quad(1\le j\le M).
\end{equation}
All shifts $b_j-y_1$ lie in $[0,\diam(Y)]$. In each term of the determinant expansion, write $n_i=e_i-j_i$, where the $j_i$ are a permutation of $0,\ldots,M$. Then
\[
 \sum_i(n_i+1)=N+\sum_i e_i-\sum_{j=0}^M j=2N-1 .
\]
It follows that $t^{(2N-1)/2}$ times each determinant equals the determinant of its entrywise normalized factors $t^{(e_i-j+1)/2}F_{e_i-j}$. Lemma~\ref{lem:llt}, applied at $u_i-(b_j-y_1)/\sqrt t$, and uniform continuity of the finitely many $\Phi_{j-e_i}$ on compact intervals give the same limit $Q_N(u)$ for every summand, with a uniform additive error when $\diam(Y)/\sqrt t\le\eta(t)$. Dividing by $\gamma_N(Y)$ averages these errors rather than multiplying them by the possibly growing number of chains. This proves \eqref{eq:joint-llt}.

For summation over the entire chamber, take $\diam(Y)/\sqrt t\le1$. In every product, the first row has an exponentially decaying left tail for every $n_1$; the last row has $n_N=M-j_N\ge0$ and hence an exponentially decaying right tail. The remaining factors and the opposite tails are bounded by polynomials. By \eqref{eq:lltbound}, the normalized determinants, and hence their average, are bounded by
\[
 C(1+\|u\|_\infty)^C \exp\{-u_1^- -u_N^+\}.
\]
On the ordered chamber $u_1^-+u_N^+\ge\|u\|_\infty$. This majorant is integrable; its lattice tails, with cell volume $t^{-N/2}$, tend uniformly to zero as the cutoff tends to infinity (compare each cell to a bounded translate of the same polynomial--exponential envelope). Diagonal faces have zero Lebesgue measure. Compact Riemann sums and the tail bound now give, for each bounded continuous $f$,
\[
 \frac{t^{M/2}}{\gamma_N(Y)}\sum_{X} f(u)K_t(Y,X)
 \ \longrightarrow\ \int_{u_1<\cdots<u_N} f(u)Q_N(u)\,du .
\]
The integrand $Q_N$ is continuous and nonnegative, by \eqref{eq:joint-llt} and approximation by ordered lattice endpoints. Taking $f=1$ for one fixed $Y$ and comparing with Theorem~\ref{thm:decay}(c) identifies its integral as $\kappa_N>0$. Taking $f$ to depend only on $u_N$ and using part~(d) of the same theorem identifies the last marginal (the fixed centring shift $y_1/\sqrt t$ vanishes). The estimates just proved give \eqref{eq:joint-survival} also for moving $Y$; division by it proves both conditional assertions.
\end{proof}

\subsection{Proof of Theorem~\ref{thm:overtaking}}\label{subsec:overtaking-proof}
\begin{proof}
Before $\tau_m$, write $L_1<\cdots<L_m<Z$ for the retained second-class ranks and the leftmost first-class position. Proposition~\ref{prop:projB} treats one first-class particle without killing. To compare the killed generators for any number of first-class particles, use the clocks at the current ranks, not at fixed labels. 
\begin{center}\small
\begin{tabular}{@{}p{.17\linewidth}@{\hspace{.025\linewidth}}p{.76\linewidth}@{}}
\toprule
Clock & Effect on the retained state \\ \midrule
$L_i$, right &
Let $b$ be the first site strictly to the right of $L_i$ outside $\{L_1,\ldots,L_m\}$; thus $b\le Z$. If $b<Z$, replace $L_i$ by $b$ and reorder. If $b$ is vacant the particle lands there. If $b$ contains an unretained second-class particle, then $R\ge m+1$ before the jump: even a landing beyond $Z$ leaves at least $m$ behind it. If $b=Z$, then $R=m$ and this jump kills.\\[3pt]
$L_i$, left &
Move to $L_i-1$ unless it is the preceding retained site, in which case the attempt is blocked.\\[3pt]
$Z$, right &
The leftmost first-class position becomes $Z+1$: if that site was already first-class it supplies the new leftmost rank; otherwise the moving first-class particle occupies it.\\[3pt]
$Z$, left &
Kill if $Z-1=L_m$. Otherwise $Z$ moves to $Z-1$ and the retained lower ranks do not change, including when an unretained second-class particle is exchanged.\\[3pt]
Other clocks &
Unretained second-class particles cannot enter the retained ranks from the right. A first-class particle to the right of $Z$ cannot move to or past $Z$ from the right. None of these clocks kills or changes the retained state.\\
\bottomrule
\end{tabular}
\end{center}
The surviving updates and the rate into the cemetery state are therefore those of the $(m+1)$-particle system killed at its first exchange. The rates are bounded, so equality of these generators identifies the laws. Proposition~\ref{prop:killed-identification} gives (a); Theorems~\ref{thm:sN} and~\ref{thm:decay} give (b) and (c).
\end{proof}

\section{Conclusions and open questions}\label{sec:open}
For fixed $N$, the same positive count $\gamma_N(Y)$ normalizes the joint local limit for $\diam(Y)=o(\sqrt t)$ and gives the Doob transform of the killed chain. For fixed initial sites, the physical and auxiliary conditional joint limits share their last marginal. They coincide for $N=2$ and differ for every $N\ge3$.

Two regimes remain to be determined. First, for $N\ge3$, what replaces \eqref{eq:joint-survival} when some initial gaps are of order $\sqrt t$ and others remain fixed?

Second, consider two overtaking thresholds $R(0)\ge m>j\ge1$ and times $0<s<t$. Let $K_s^{[m+1]}$ denote the killed kernel of Proposition~\ref{prop:killed-identification} with $m+1$ particles. The Markov property at time $s$ and Theorem~\ref{thm:overtaking} give
\[
 \PP(\tau_m>s,\tau_j>t)=
 \sum_{z_1<\cdots<z_{m+1}} K_s^{[m+1]}(Y^{(m)},z)\,
 \mathcal S_{j+1}\big(t-s;(z_1,\ldots,z_j,z_{m+1})\big).
\]
On $\{\tau_m>s\}$ the first $j$ retained positions and the last one determine the future $j$-threshold survival probability. Can this sum be reduced to a single determinant, or evaluated asymptotically as $t\to\infty$ with $s/t\to\theta\in(0,1)$?


\section*{Use of AI tools}
In preparing this work, the authors utilized ChatGPT Astra (OpenAI) and Claude Fable 5.2 (Anthropic) for proof verification, literature search, Python code development and validation, figure and table generation, as well as text revision and editing. After using these tools, the authors thoroughly reviewed and edited all generated outputs and assume full responsibility for the mathematical arguments, computations, references, and final text.

\section*{Funding}
Eunghyun Lee was supported by Nazarbayev University under the Faculty-Development Competitive Research Grants Program for 2024–2026 (Grant No. 201223FD8822, EL). Yuncherl Choi was supported by the National Research Foundation of Korea (NRF) grant funded by the Korea government (MSIT) (No. RS-2021-NR063084).

\appendix
\section{Proof of Theorem~\ref{thm:matrix}}\label{app:matrix}
\begin{proof}
Use a reduced chain adapted to the distinguished letter, as in \cite{Lee-Raimbekov-2022}.

(a) Let $\rho=\rho_1\cdots\rho_{N-1}:=\sigma\setm1$, a word in the letters $2,\dots,N$. \emph{Stage 1.} Starting from the identity, we rearrange the letters $2,\dots,N$, which occupy the slots $2,\dots,N$, into $\rho_1\cdots\rho_{N-1}$, by a chain of $|\Inv(\rho)|$ elementary steps at slots $\ge2$. The pairs exchanged are the elements of $\Inv(\rho)$, and the result is the word $1\rho_1\cdots\rho_{N-1}$. \emph{Stage 2.} We apply the elementary steps at slots $1,2,\dots,r-1$, in this order; the step at slot $i$ exchanges the pair $(\rho_i,1)$. The result is $\rho_1\cdots\rho_{r-1}\,1\,\rho_r\cdots\rho_{N-1}=\sigma$. The total number of steps is $|\Inv(\rho)|+r-1=|\Inv(\sigma)|$, so that
\[
\bA_\sigma=\bT_{r-1}(\rho_{r-1},1)\cdots\bT_1(\rho_1,1)\,\mathbf B,
\]
where $\mathbf B$ is the product of the operators of Stage~1. Each of these acts on two slots in which $\wa_1=12\cdots2$ has the letters $22$; by \eqref{eq:Taction}, $\mathbf Be_{\wa_1}=S_\rho\,e_{\wa_1}$. We claim that for $0\le k\le r-1$,
\begin{equation}\label{eq:indA}
\bT_k(\rho_k,1)\cdots\bT_1(\rho_1,1)\,e_{\wa_1}=\prod_{i=1}^kT_{\rho_i1}\;e_{\wa_{k+1}}-\sum_{m=1}^k\ \prod_{i=1}^{m-1}T_{\rho_i1}\prod_{i=m+1}^{k}S_{\rho_i1}\;e_{\wa_m}.
\end{equation}
For $k=0$ both sides equal $e_{\wa_1}$. Assume \eqref{eq:indA} for some $k<r-1$ and apply $\bT_{k+1}(\rho_{k+1},1)$, which acts on the slots $k+1,k+2$. The word $\wa_{k+1}$ has the letters $12$ in these slots, so $e_{\wa_{k+1}}\mapsto-e_{\wa_{k+1}}+T_{\rho_{k+1}1}e_{\wa_{k+2}}$; for $m\le k$ the word $\wa_m$ has the letters $22$ there, so $e_{\wa_m}\mapsto S_{\rho_{k+1}1}e_{\wa_m}$. This gives \eqref{eq:indA} with $k+1$ in place of $k$. Taking $k=r-1$ and recalling that $\rho_i=\sigma(i)$ for $i<r$, we obtain \eqref{eq:matA}.

(b) Let $\rho=\rho_1\cdots\rho_{N-1}:=\sigma\setm N$, a word in the letters $1,\dots,N-1$. \emph{Stage 1.} By $|\Inv(\rho)|$ elementary steps at slots $\le N-2$ we transform the identity into $\rho_1\cdots\rho_{N-1}N$. \emph{Stage 2.} We apply the elementary steps at slots $N-1,N-2,\dots,r$, in this order; the step at slot $i$ exchanges the pair $(N,\rho_i)$. The result is $\rho_1\cdots\rho_{r-1}\,N\,\rho_r\cdots\rho_{N-1}=\sigma$, so that $\rho_i=\sigma(i+1)$ for $i\ge r$, and
\[
\bA_\sigma=\bT_{r}(N,\rho_{r})\cdots\bT_{N-1}(N,\rho_{N-1})\,\mathbf B .
\]
The operators of Stage~1 act on slots in which $\wb_N=1\cdots12$ has the letters $11$, hence $\mathbf Be_{\wb_N}=S_\rho\,e_{\wb_N}$. The operator $\bT_i(N,\rho_i)$ acts on the slots $i,i+1$, where $\wb_{i+1}$ has the letters $12$ and $\wb_m$, $m>i+1$, has the letters $11$. At $k=N$ the empty product gives $e_{\wb_N}$. For the step from $k+1$ to $k$, \eqref{eq:Taction} sends $e_{\wb_{k+1}}$ to $-e_{\wb_{k+1}}+T_{N\rho_k}e_{\wb_k}$ and each $e_{\wb_m}$ with $m>k+1$ to $S_{N\rho_k}e_{\wb_m}$. Substitution in the formula at $k+1$ gives the following identity for $r\le k\le N$:
\[
\bT_k(N,\rho_k)\cdots\bT_{N-1}(N,\rho_{N-1})\,e_{\wb_N}=\prod_{i=k}^{N-1}T_{N\rho_i}\;e_{\wb_k}-\sum_{m=k+1}^{N}\ \prod_{i=m}^{N-1}T_{N\rho_i}\prod_{i=k}^{m-2}S_{N\rho_i}\;e_{\wb_m},
\]
and \eqref{eq:matB} follows by taking $k=r$ and substituting $\rho_i=\sigma(i+1)$.
\end{proof}

\section{Proofs of Lemmas~\ref{lem:sum} and~\ref{lem:llt}}\label{app:lemmas}
\begin{proof}[Proof of Lemma~\ref{lem:sum}]
Absolute convergence follows from $|\nabla^iG_j(x)|\le 2^iC_b R^{x}$ for $x\le b$, enlarging $C_b$ over the finitely many positive sites if necessary. We argue by induction on $n$; for $n=1$ the left-hand side telescopes to $G_j(b)$. Let $n\ge2$ and write $g=g_J$. Expanding the determinant along its last row and applying the induction hypothesis (with $b$ replaced by $x_n-1$) to each cofactor, the sum over $x_1<\dots<x_{n-1}<x_n$ equals $\det\big(g(x_n-n+1),\dots,g(x_n-1),\nabla^ng(x_n)\big)$. By the binomial formula, $\nabla^ng(x)-g(x)-(-1)^ng(x-n)$ is a linear combination of $g(x-1),\dots,g(x-n+1)$. Hence, with $V(x):=\det\big(g(x-n+1),\dots,g(x)\big)$,
\[
\det\big(g(x-n+1),\dots,g(x-1),\nabla^ng(x)\big)=V(x)-V(x-1),
\]
where we used that moving the row $g(x-n)$ from the last to the first place produces the sign $(-1)^{n-1}$. Since $V(x)\to0$ as $x\to-\infty$, the sum over $x_n\le b$ is $V(b)$.
\end{proof}

\begin{proof}[Proof of Lemma~\ref{lem:llt}]
The estimates below hold for all $p\in[0,1]$, with $q=1-p$.
\emph{The case $n\ge0$.} The integrand in \eqref{eq:Fdef} is then analytic on $\C\setm\{0\}$, so the contour may be any circle $|\xi|=\rho$. Taking $\rho=1$ and $\xi=e^{-i\theta}$ we get
\[
F_n(x;t)=\frac1{2\pi}\int_{-\pi}^{\pi}(1-e^{i\theta})^n\,e^{-i\theta x+t\psi(\theta)}\,d\theta,\qquad \psi(\theta):=pe^{i\theta}+qe^{-i\theta}-1 ,
\]
and, substituting $\theta=r/\sqrt t$,
\[
t^{(n+1)/2}F_n(x;t)=\frac1{2\pi}\int_{-\pi\sqrt t}^{\pi\sqrt t}\big[\sqrt t\,(1-e^{ir/\sqrt t})\big]^n\,e^{-iru}\,e^{t(\psi(r/\sqrt t)-ivr/\sqrt t)}\,dr .
\]
As $t\to\infty$ the integrand without the factor $e^{-iru}$ converges pointwise to $(-ir)^ne^{-r^2/2}$, and its modulus is at most $|r|^ne^{-2r^2/\pi^2}$, because $|1-e^{i\theta}|\le|\theta|$ and $\operatorname{Re}\psi(\theta)=\cos\theta-1\le-2\theta^2/\pi^2$ for $|\theta|\le\pi$. By dominated convergence, $t^{(n+1)/2}F_n(x;t)-\frac1{2\pi}\int_{\mathbb R}(-ir)^ne^{-iru-r^2/2}dr\to0$ uniformly in $x$, and the last integral is $\varphi^{(n)}(u)=\Phi_{-n}(u)$. This proves \eqref{eq:lltlimit}. For \eqref{eq:lltbound} we take $\rho=e^{-b}$ with $b:=\operatorname{sgn}(u)/\sqrt t$ (and $b:=1/\sqrt t$ if $u=0$), so that $|b|\le1$. For $\xi=\rho e^{i\theta}$ we have $\operatorname{Re}\eps(\xi)=\lambda(\cos\theta-1)+\lambda-1$ with $\lambda:=pe^{b}+qe^{-b}$, and $\lambda\ge e^{-1}$, $\lambda-1\le vb+b^2$ (use $e^s\le1+s+s^2$ for $|s|\le1$). Hence
\[
|\xi|^{x}e^{t\operatorname{Re}\eps(\xi)}\le e^{-b(x-vt)+tb^2}\,e^{-t(1-\cos\theta)/e}\le e^{1-|u|}\,e^{-2t\theta^2/(e\pi^2)} .
\]
Moreover $|h(\xi)|=|1-e^{b}e^{-i\theta}|\le e\,(t^{-1/2}+|\theta|)$. Therefore
\begin{align*}
|F_n(x;t)|&\le\frac1{2\pi}\int_{-\pi}^{\pi}|h(\xi)|^n|\xi|^xe^{t\operatorname{Re}\eps(\xi)}d\theta\\
&\le\frac{e^{n+1}}{2\pi}\,e^{-|u|}\int_{-\pi}^{\pi}(t^{-1/2}+|\theta|)^n\,e^{-2t\theta^2/(e\pi^2)}d\theta\ \le\ C_n\,t^{-(n+1)/2}e^{-|u|}.
\end{align*}

\emph{The case $n\le-1$.} We use downward induction on $n$, assuming the lemma for $n+1$. By Lemma~\ref{lem:F}(e), $F_n(x;t)=\sum_{x'\le x}F_{n+1}(x';t)$, that is,
\[
t^{(n+1)/2}F_n(x;t)=\int_{-\infty}^{u}f_t(r)\,dr,\qquad f_t(r):=t^{(n+2)/2}F_{n+1}\big(\lceil vt+r\sqrt t\,\rceil;t\big).
\]
By the induction hypothesis $|f_t(r)|\le C\beta_{n+1}(r)$ (note that $\beta_{n+1}(r')\le C'\beta_{n+1}(r)$ if $|r'-r|\le1$), and $\int_{-\infty}^u\beta_{n+1}(r)dr\le2\beta_n(u)$; this gives \eqref{eq:lltbound}. Also by the induction hypothesis, $f_t\to\Phi_{-n-1}$ uniformly on compact sets, which together with the bound $|f_t|\le C\beta_{n+1}$ gives $\int_{-\infty}^uf_t(r)dr\to\int_{-\infty}^u\Phi_{-n-1}(r)dr=\Phi_{-n}(u)$ uniformly for $|u|\le A$.
\end{proof}

\section{Numerical checks}\label{app:num}
For $Y=(0,2,5)$, $p=0,0.6,1$ and $t=1,2$, the finite-window killed generator agrees with the determinant at the selected endpoints and prey marginals to within $2\times10^{-15}$. Other checks include 2,184 projected-rate comparisons, 4,280 rational diagonal-coefficient tests and 70 exact tests of \eqref{eq:joint-chains}. For the two crossing marginals \eqref{eq:FN} and \eqref{eq:GN}, all 2,970 Taylor-coefficient comparisons through order eight agree with direct iteration of the two-species generator. The Vandermonde difference in Remark~\ref{rem:initial-differences} agrees with the chain count on all $492$ configurations of $2\le N\le7$ sites in $\{0,\ldots,8\}$.

\raggedbottom
\bibliographystyle{abbrv}
\bibliography{refs}

\end{document}